\documentclass[a4paper,11pt,]{article}
\usepackage{amssymb,amsmath}
\usepackage{natbib}
\newtheorem{defin}{Definition}
\newtheorem{defin-theor}{Definition-Theorem}
\newtheorem{defin-rem}{Definition-Remark}
\newtheorem{rem}{Remark}
\newcommand{\rf}[1]{(\ref{#1})}

\newcommand{\bc}{\begin{center}}
\newcommand{\ec}{\end{center}}
\newcommand{\be}{\begin{equation}}
\newcommand{\ee}{\end{equation}}
\newcommand{\bea}{\begin{eqnarray}}
\newcommand{\eea}{\end{eqnarray}}

\newcommand{\bfr}{\begin{flushright}}
\newcommand{\efr}{\end{flushright}}
\newcommand{\bfl}{\begin{flushleft}}
\newcommand{\efl}{\end{flushleft}}
\newcommand{\defeq}{\stackrel{\mbox{\scriptsize def}}{=}}
\newcommand{\comp}{\!\circ\!}
\newtheorem{theor}{Theorem}

\begin{document}
\title{Dynamics of Birational Plane Mappings.\\ The Arnold complexity difference equation.}
\author{Konstantin V.~Rerikh\footnote{e-mail:rerikh@thsun1.jinr.ru}}
\date{}
\maketitle
\begin{center}
{Bogoliubov Laboratory of Theoretical Physics, JINR,\\
141980, Dubna, The Moscow Region, Russian Federation\\}
\end{center}

\begin{abstract}
{\rm
We consider a dynamics of a generic birational plane map $\Phi_n: {\bf\rm CP}^2\to{\bf\rm CP}^2$,
${\bf\rm CP}^2$ -image  of the birational mapping (inverse map is also rational)
 $F_n : {\bf\rm C}^2\to {\bf\rm C}^2$
and its such important characteristic as the Arnold complexity $C_A(k)$, which is proportional
$d(k)={\rm deg}(\Phi_n^k)$- a degree of $k-$iteration of the map $\Phi_n$, on the basis
 on algebraic-geometrical properties of such maps. Additional importance of this characteristic
 follows from  the Veselov conjecture about the polynomial boundedness of the growth of $d(k)$ for
 integrable dynamical systems with a discrete time defined by birational plane maps.
  The autonomous linear difference equation with integer coefficients
 for $d(k)$ is obtained. This equation is fully defined by $\sigma_1$ nonnegative integers
 $m_1,\cdots, m_{\sigma_1}$ that are  determined  by relations:
  $\Phi_n^{-m_i}(O_{\alpha_i})=O^{(-1)}_{\beta_i}, \quad i\in(1,2,\cdots,\sigma_1)$,
  where $\Phi_n^{-m_i}$ is
  $m_i$-iteration of inverse map, $O_{\alpha_i}, O^{(-1)}_{\beta_i},\quad
  \alpha_i, \beta_i \in (1,2,\cdots,\sigma)$ are indeterminacy points
  of the direct and inverse maps,  $\sigma_1\leq\sigma$ and $\sigma$ is a number of indeterminacy
   points of $\Phi_n,\Phi_n^{-1}$. If $\sigma_1$ is equal to zero that $d(k)=n^k$, otherwise
   the growth of $d(k)$ is fully defined by a root spectrum of the secular equation associated
    with the difference equation for $d(k)$. The Veselov conjecture corresponds to
    the root spectrum consisting of values being equal to modulo one.
    The author doesn't suppose that the reader has acquaintance with the algebraic
geometry (AG) in ${\bf\rm CP}^2$ and the dynamical systems theory (DST) or the functional
equations since in the paper there are given all needed definitions of used concepts of AG or
 DST and theorems.
}
\end{abstract}


 \section{Introduction. Set of the problem and Main Result.}

Let consider the system of birational functional equations (BFEs) for
 functions $y(w):~{\rm\bf C}\to{\rm\bf C}^{N} $ in one complex variable $w$
  of the form
\be y(w+1)=F_n(y(w)),\quad y(w):~{\rm\bf C}\to{\rm\bf C}^{N},
 \quad w\in {\rm\bf C}, \quad F_n \in {\rm\bf Bir}({\rm\bf C}^{N}).
\label{eq:bfes1}
 \ee
  For $w=m\in{\rm\bf Z}$ the above BFEs  are a dynamical system
 with discrete time or cascade.
 Here the map
$
  F_n:~y\mapsto y'= F_n(y)=\frac{f_i(y)}{f_{N+1}(y)}, i=(1,2,\ldots,N)$,
$f_i(y)~ \mbox{for}~ \forall~ i$ are polynomials in $y$,
  ${\rm deg}F_n(y)= \max_{i=1}^{N+1} \left\{ {\rm deg}(f_i(y)) \right\}=n,
$
is a given birational map of the group of all automorphisms of ${\rm\bf C}^N \to
 {\rm\bf C}^N$\\ ( the Cremona group or $ {\rm\bf Bir}({\rm\bf C}^{N}$))
  with the coefficients from $ {\bf\rm C}$.

A preliminary investigation of the dynamics of the mapping $F_n$ is important in the context of
consideration of the integrability problem of  dynamical systems or the BFEs of the form
\rf{eq:bfes1}. Such consideration is more  convenient and effective to realize in
${\rm\bf CP}^{N}$. Let us give the definition of the mapping $\Phi_n$, the  image of the map $F_n$
in ${\rm\bf CP}^{N}$, at $N=2$ since below  we shall mainly consider the dynamics of plane birational
mappings.
At the transition to ${\rm\bf CP}^2$ $y\mapsto z: y_i = z_i/z_3, i=(1, 2)$ the
maps $F_n, F_n^{-1}$ transform into
the maps $\Phi_n, \Phi_n^{-1} $:
\bea
 \Phi_n:&~& z\mapsto z', \quad z_1':z_2':z_3'= \phi_1(z): \phi_2(z):\phi_3(z),\quad z, z'\in
{\rm\bf CP}^2,
 \label{eq:fin2} \\
\phi_i(z) & = & z_3^nf_i(z_l/z_3), \quad i\in(1, 2, 3),\quad l=(1, 2),
\eea
 and $\phi_i(z)$ are
 homogeneous polynomials in $z$ without any common factors. The map
 \be
\Phi_n^{-1}: z'\mapsto z,~ z_1:z_2:z_3= \phi^{(-1)}_1(z'):
\phi^{(-1)}_2(z'):\phi^{(-1)}_3(z'),~ z, z'\in {\rm\bf CP}^2
 \label{eq:fin2'}
 \ee
  is defined analogously.

In the abstract and above we used such familiar concepts  as "integrability", "integrable
maps", "integrable dynamical systems", and "integrable functional equations". In order  to
avoid different understanding of these terms, we shall below give our definition of these concepts.
The integrability problem for dynamical systems and functional equations is  solved
if we obtain for the BFEs \rf{eq:bfes1} the family of first integrals $I(y(w),w)$ of dimension
of $1\leq m \leq N-1$ of the form:
\bea
I(y(w),w) & = & c(w), c(w):{\bf\rm C}\to {\bf\rm C}^{m},\quad c(w+1)=c(w), \label{eq:c(w)}\\
I(y,w) & : & ({\bf\rm C}^N\otimes{\bf\rm C}) \to {\bf\rm C}^m, \label{eq:I}
\eea
where arbitrary periodic functions $c(w)$ in $w$ parameterize the level lines of first integrals.

In addition to these $m$ first integrals we always have one more first transitive integral parameterized
by the periodic function
\be
w\to w+\beta(w),\quad \beta(w):{\bf\rm C}\to {\bf\rm C},\quad \beta(w+1)=\beta(w). \label{eq:beta(w)}
\ee
If integer $m=N-1$, we can speak about full integrability, otherwise partial integrability of equation
 \rf{eq:bfes1} at $1 \leq m<N-1$. We have a general solution of equations BFEs \rf{eq:bfes1}
 if we obtain the solution of BFEs \rf{eq:bfes1} in the explicit form $y=Y(w,c(w))$ with $c(w):
 {\bf\rm C}\to{\bf\rm C}^{N-1}$. If $I(y,w)$ is a rational function of $(y, w)$ or a rational
 function of $y$ and fraction-linear one in variables $\tau(w)$ where $\tau(w)$ are variables of
  the form
 $\lambda^w $, then we can speak about the algebraic integrability of BFEs \rf{eq:bfes1}.
 If $I(y,w)\in
 Hol(y,w)$, i.e. $I(y,w)$ is a holomorphic function of variables, then we can speak about non-algebraic
 integrability of BFEs \rf{eq:bfes1}. If the algebraic integrability is the subject for using
 algebra-geometrical methods, then non-algebraic integrability is the subject for using classical
 methods and theorems from the theory of dynamical systems  due to H. Poincare, C.L.~Siegel,
 G.D.~Birkhoff, A.N.~Kolmogorov, V.I.~Arnold, J.~Moser, D.V.~Anosov and others
  (see \cite{DS-1-88a}, \citet{DS-1-88b}), and also using  classical results of the number theory and
  the transcendental number theory (see, for example, \citep{moser-94}, \citep{rerikh-95a},
  \citep{rerikh-97}, \citep{rerikh-98b} as examples of using classical results of A.~Baker
  \citep{baker-90}, \citep{baker-66},  \citep{baker-67a,baker-67b,baker-68}, \citep{baker-71},
 \citep{baker.wustholz-93} and N.I.~Feldman
   \citep{feldman-82-rus}, \citep{feldman-68-rus}).  In the paper,
 we shall not  discuss the integrability problem in more detail.
There are also other definitions of the concept integrability. The A.P.~Veselov definition as
 applied to the dynamical systems of the form \rf{eq:bfes1} acting in the plane (bipolynomial
 maps, \citep{veselov-89}, \citep{veselov-91}) is as follows: "The map $F_n$ is integrable if
  there exist another map $\Psi$ for which $\Psi^m \neq F_n^k\forall k,m \in{\bf\rm Z}$,
   where $\Psi^m$ is the $m$th iteration of the map $\Psi$ but $F_n^k$ is the $k$th iteration
  of the map $F_n$, commuting with $F_n$: $\Psi\comp F_n\equiv F_n\comp\Psi$ ".
  The Moser definition means the existence of a holomorphic map $H(y):u=H(y),\quad {\bf\rm C}^N\to
  {\bf\rm C}^N$ that transforms the maps $F_n$ to its linear part $H\comp F_n\comp H^{-1}\equiv
  A\defeq\frac{\partial F_n}{\partial y}|_{y=y_0}$, where $y_0$ is a fixed point of the map $F_n$
  but the matrix $A$ defines the linear part of the map $F_n$. Such a definition is natural for
  the theory of dynamical systems. (see examples of non-algebraically integrable dynamical systems
  \citep{rerikh-92}, \citep{rerikh-95b} \citep{rerikh-95a}, \citep{rerikh-97}, \citep{rerikh-98a}). This Moser definition
  of integrability is in fact a local concept in the neighbourhood of a fixed point of a map as
  well as a concept of a local non-integrability in a neighbourhood of a fixed point of a map.
  (See  the Moser example \citep{moser-60} of non-integrable cubic bipolynomial map in the neighbourhood of the zero
  elliptic fixed point.)

V.I. Arnold in papers \citep{arnold-90a}, \citep{arnold-90b}  introduced
 and  investigated  such a characteristic of a dynamical system as the topological
{\it complexity}~ of the intersection of a submanifold $X$ of manifold $M$,
 moved by a dynamical system, a smooth mapping $A:~~M\to M$ with the other
 given compact smooth submanifold $Y$ of $ M$:
$$
Z_k=(A^k X)\cup Y.
$$

    In the simplest case for
 plane mappings $\Phi$ the complexity $C_A^{\Phi}(k)\equiv Z_k$
  can be defined
\citep{veselov-92} as the  number of intersection points of a fixed
 curve $\Gamma_1$ with
 the image of another curve $\Gamma_2$ under the $k^{th}$ iteration of
$\Phi$:
$$
 C^{\Phi}_{A;\Gamma_1\Gamma_2}(k) = \#\bigl(\Gamma_1 \cap \Phi^k(\Gamma_2)\bigr).
$$
If the mapping $\Phi$ is a birational one from the group ${\rm \bf BirCP}^2$ and the curves
$\Gamma_1,$ $\Gamma_2$ are algebraic curves in ${\rm \bf CP}^2$, then it is easy
 to see that the growth of $ C^{\Phi}_{A;\Gamma_1\Gamma_2}(k)$ will in general
be as follows:
$$
 C^{\Phi}_{A;\Gamma_1\Gamma_2}(k) = {\rm deg}(\Gamma_1) {\rm deg}(\Gamma_2)
 d_{\Phi}(k)\leq {\rm deg}(\Gamma_1) {\rm deg}(\Gamma_2) ({\rm deg}\Phi)^k,
$$
where $d_{\Phi}(k)={\rm deg}(\Phi^k)$ is the degree of the mapping $\Phi^k
=\Phi\comp\Phi \comp\cdots \comp \Phi$,
 which agrees well
 with general Arnold's results for smooth mappings
 and diffeomorphisms \citep{arnold-90a}, \citep{arnold-90b}.

 The Arnold complexity was found to be an important characteristic in the  context  of
  the integrability of such dynamical systems. A.P.~Veselov introduced in \citep{veselov-91},
  \citep{veselov-92} the conjecture about
  a polynomial growth  of the Arnold complexity $d(k)$ with $k$ for integrable plane birational
   mappings and proved it for integrable bipolynomial ones in \citep{veselov-89}
   ( $d(k)$ is bounded by a constant).
    To be more exact we reformulate the A.P.~Veselov conjecture as  " all integrable birational
    mappings have a polynomially bounded growth of the Arnold complexity $d(k)$ on $k$".


 The validity of the A.P.~Veselov conjecture was also confirmed  for
   many concrete integrable mappings by different researchers   so that it is actual to prove
   it in a general case. This paper is the first step in this direction.
   The aim of this paper is to discuss the dynamics of generic
 birational plane  mappings in the frames of their algebraic-geometrical properties and
  obtain the autonomous linear difference equation for the Arnold complexity $d(k)$.

 The main result of the paper is  the obtained autonomous linear difference equation
 for $d(k)$. This equation is fully defined by $\sigma_1$ nonnegative integers
 $m_1,\cdots, m_{\sigma_1}$ that are  determined  by relations:
  $\Phi_n^{-m_i}(O_{\alpha_i})=O^{(-1)}_{\beta_i}, \quad i\in(1,2,\cdots,\sigma_1)$,
  where $\Phi_n^{-m_i}$ is the
  $m_i$-iteration of the inverse map, $O_{\alpha_i},~ O^{(-1)}_{\beta_i}, \\
  \alpha_i, \beta_i \in (1,2,\cdots,\sigma)$, are indeterminacy points
  of the direct and inverse maps,  $\sigma_1\leq\sigma$ and $\sigma$ is a number of indeterminacy
   points of $\Phi_n,\Phi_n^{-1}$. If $\sigma_1$ is equal to zero, then $d(k)=n^k$, otherwise
   the growth of $d(k)$ is fully defined by a root spectrum of the secular equation associated
    with the difference equation for $d(k)$. The A.P.~Veselov conjecture corresponds to
    the root spectrum consisting of values being equal to modulo one. Thus, this equation
 gives the possibility to present  all sets of numbers $m_1,\cdots, m_{\sigma_1}$ corresponding
 to integrable mappings if the A.P.~Veselov  conjecture is true.

In the following section,  we shall perform a brief excursus into the theory of the Cremona
transformations in the plane following  \citep{hudson-27}, \citet{topics-70},  \citep{iskovskikh.reid-91},
 \citep{shafarevich-77}, \citep{coble-61}.  In Section
\ref{decomp} we introduce a  new notion --the decomposition of  sets of the indeterminacy points
of direct $\Phi_n$ and inverse $\Phi^{-1}_n$ maps. Then in Section \ref{dynamics.diff.eq} we
obtain the main equations of the dynamics of a generic birational map and the difference equation for
the Arnold complexity d(k). Different examples for illustration of Sections \ref{alg.geom}-
\ref{dynamics.diff.eq} are set in Appendices A, B.

\section{ Brief excursus into the algebraic geometry} \label{alg.geom}

 Let $z=(z_1,z_2,z_3)$ be a point of the projective plane ${\bf\rm CP^2}$.
Let us consider a general curve of degree $\mu$ $f_\mu(z)$ defined by
 the equation
\be
f_\mu(z)=\sum_{|l|=\mu}c_lz^l=0,\quad l\defeq (l_1,l_2,l_3),\quad |l| \defeq l_1+l_2+l_3,
\label{eq:fa1}
\ee
which has, in general, $\frac{\mu(\mu+3)}{2}$ free parameters.

\subsection{\bf Linear systems of curves}\label{LSC}

\begin{defin} \label{mult.point}
{\rm Let $P=(z^*_1, z^*_2, z^*_3)$ be a point of the curve
\rf{eq:fa1} and let $z_3^*$ be the coordinate of the point $P$
which is nonzero, but, therefore, we can assign $z^*_3=1$ as
a result of the change $P \mapsto P/z^*_3 $. The point $P$ is called
an $r$-fold one of the curve \rf{eq:fa1} if $f_\mu(z)$ has the
following form in the system of coordinates
$z':~~z'_1=z_1-z^*_1z_3, z'_2=z_2-z^*_2z_3, z'_3=z_3 $
\bea
f_\mu(z) & = & f'_\mu(z')=\sum_{k=r}^{\mu}{z'_3}^{\mu-k}u_k(z'_1, z'_2),
\label{eq:fa2} \\
r & \defeq & \mbox{mult}(f_\mu(z))|_{z=z^*},\quad \mbox{mult}\defeq \mbox{multiplicity},
\eea
 where $u_k(z')$ are homogeneous polynomials of
degree $k$ in variables
 $z'_1, z'_2$, but the first function $u_r(z')$ in expansion \rf{eq:fa2}
 defines $r$ tangents for the curve at the point $P$.

An $r$-fold point imposes $\frac{r(r+1)}{2} =\sum_{i=0}^{r-1}(i+1)$
conditions ensured  for a curve \rf{eq:fa1} of the form \rf{eq:fa2} and is called a simple, double, triple one, if $r$
is equal $1, 2, 3$ and so on.  $ \triangleleft$
}
\end{defin}

 Let us give a definition of a linear system of curves which is important
 in what follows.

\begin{defin} \label{linear.sys}
 {\rm
(For more details see \citet{topics-70} and also all references therein
 on
 results and notions reviewed here.)~~~ The system of plane curves
 ${f_\mu}$ of degree $\mu$ is represented by an equation of the form
\be
f_\mu=\sum_{i=1}^{k+1}c_if_i(z)= 0,
\label{eq:fa3}
\ee
where  the functions $f_i(z)$ are homogeneous polynomials from
 $z=(z_1, z_2, z_3)$ of the same order $\mu$ and are linearly independent, is
a linear system of curves (LSC) of dimension $k$.

\begin{defin}\label{base.LSC}
{\rm
 A point $B^{(r_j)}_j$ which is at least an $r_i$-fold one for each curve of the system is called
an $r_i$-fold basis point (see Definition \ref{mult.point}) but a join of all basis points is
called a basis set or a base of the LSAC $ {\bf B}=\bigcup_j B^{(r_j)}_j $, where $
j\in{1,2,\cdots, N_B} $. Thus, the base is fully defined by  two sets:  the set of basis points
and the set of theirs multiplicities on the LSAC which are linked  with each other, as it is set
above.}
 \end{defin}

 Counting separately the
conditions imposed by all basis points necessary for reduction from
the general curve \rf{eq:fa1} to the LSAC \rf{eq:fa3}, we have the virtual dimension
$$
K=\frac{\mu(\mu+3)}{2}-\sum_{j=1}^{N_{{\bf\rm B}}}\frac{r_j(r_j+1)}{2}.
$$
In certain cases the conditions are not independent so that the effective
dimension is
\be
k=K+s,
\label{eq:fa4}
\ee
where $s$ is the number of independent relations among the linear conditions imposed
 by the base $B$ on curves of order $\mu$. A system for which $s=0$ is said
 to be regular, otherwise irregular with irregularity \\( superabundance) $s$.

The effective genus $p$ of a general curve of the irreducible system coincides
with virtual $P$ and is
\be
p=P=\frac{(\mu-1)(\mu-2)}{2}-\sum_{j=1}^{N_{{\bf\rm B}}}\frac{r_j(r_j-1)}{2}.
\label{eq:fa5}
\ee
 For reducible curves the effective genus $p$ equals
 \be
 p=P+ c -1, \label{eq:genus.reduc}
 \ee
 where $c$ is  a number of components of a reducible curve.
The number of variable intersections of two curves of the system is
the grade $D$
\be
D=\mu^2-\sum_{j=1}^{N_{{\bf\rm B}}} r_j^2.
\label{eq:fa6}
\ee
The numbers $K, D$ and $P$ satisfy the relation
$$
K=D-P+1,
$$
so that
$$
k=D-p+s+1
$$
for an irreducible system, where $s$ is $\leq p$, because $D\geq (k-1)$.
The numbers $D, p$ and $k$ are invariant under birational mappings.
A linear system of dimension $k=1, 2, 3$ is called a pencil, a net and
a web, respectively (see also Remark \ref{non.ord.linear.sys}).
$\triangleleft$
}
\end{defin}

\begin{rem} \label{non.ord.linear.sys}
{\rm
 Definition \ref{linear.sys} can be extended to the case of the
linear systems of curves of which the system of the basis $r$-fold points
 includes
some non-ordinary (extra ordinary) $r$-fold ones (see Definition
 \ref{non.ord.point}
below).
}
\end{rem}

\begin{defin} \label{non.ord.point}
{\rm The $r$-fold basis point of the linear system of curves is
called a non-ordinary one if at this point the linear system of
curves satisfies some additional tangency conditions as the
existence of $r_1, 1\leq r_1\leq r,$
  common tangents ($r$-fold point of a simple contact) or the existence of some
fixed curve touching upon these common tangents and osculating
with each curve of the system ($r$-fold point of higher contact).
 Each $r$-fold non-ordinary point can
  be represented by the system of infinitely near ordinary points and be resolved
   using the technics
  of resolution of singularities of plane curves (see \citep{hudson-27},
   Chapter VII and also
below Section \ref{bir.equiv.res.sing} ). }
\end{defin}

Let us give a definition of a birational mapping
 $\Phi_n: {\bf\rm CP}^2\to {\bf\rm CP}^2 $.

\subsection {\bf Definition of Birational Map, Noether theorem and Quadratic maps}

 \begin{defin}\label{map} {Birational map.}
{\rm A mapping $\Phi_n: z\mapsto z',~ z,z'\in{\bf\rm CP}^2\ $ in
\rf{eq:fin2},
 where $\phi_i$ are homogeneous polynomials in $z, i=(1,2,3)$, of degree $n$,
is called a birational mapping if it assigns one-to-one correspondence
between $z$ and $ z' $, while the inverse mapping is given by
\rf{eq:fin2'}
and it is also rational (genus $p=0$),
$\phi_i'$ being also homogeneous polynomials in $z'$, moreover, $\phi_i$ and
 $\phi^{(-1)}_i$ have no common factors.

Associated with $\Phi_n$ and $\Phi^{-1}_n$  the linear systems of curves
 $\phi, \phi^{(-1)}$  of dimension $k=2$, genus $p=0$ and grade $D=1$
\bea
\phi = c_1\phi_1 +c_2 \phi_2 + c_3 \phi_3, \label{eq:fa7}  \\
\phi^{(-1)} = c^{(-1)}_1\phi^{(-1)}_1 +c^{(-1)}_2 \phi^{(-1)}_2 + c^{(-1)}_3 \phi^{(-1)}_3  \label{eq:fa8}
 \eea
  (for $ c_i, c^{(-1)}_i \in {\bf\rm C} $) are fully given  by theirs  bases ${\bf B, B}^{-1}$
   (for  bases
of  LSACs associated with  maps we shall use symbols ${\bf B}\defeq {\bf O}, {\bf B}^{-1}\defeq
 {\bf O}^{-1}, r_\alpha\defeq i_\alpha, r^{(-1)}_{\beta}\defeq i^{(-1)}_\beta, \alpha \in(1,2,\cdots ,\sigma\defeq
  N_{\bf B}), \beta \in
(1,2, \cdots, {\sigma}^{(-1)}\defeq N_{{\bf B}^{-1}}),
(\sigma={\sigma}^{(-1)}). $) define the first and second rational
nets which are images of nets of lines. The basis points $O_\alpha,
O^{-1}_\beta$ are indeterminacy ones for the maps $\Phi_n, \Phi^{-1}_n$
and are called fundamental ones (or F-points). The equality of genus
$p$ to zero is a necessary and sufficient condition for the
birationality of the rational map $\Phi_n$ \rf{eq:fin2}.
$\triangleleft$ }

\end{defin}

\begin{theor}  \label{Noether} (M. Noether)
{\rm Every Cremona plane mapping $\Phi_n$ \rf{eq:fin2} can be resolved into
quadratic
 mappings
$$
\Phi_n=C\comp Q_1\comp Q_2 \cdots \comp Q_j ,
$$
where $C$ is a collineation (linear mapping in ${\bf\rm CP}^2$), but mappings
 $Q_1,\cdots, Q_j $ are quadratic ones.
. $\triangleleft$ }
\end{theor}

At the end, we should give the definition of the main object-- the gene\-ra\-tors of the Cremona
group, namely, birational quadratic mappings.

\begin{defin} \label{quadratic}
{\rm Any generic quadratic Cremona mapping is generated by a composition
 \be
  \Phi_2 \equiv B^{-1}\comp I_s\comp B_1, \label{eq:quad}
 \ee
  where
   \be
   B: z\mapsto j^{(-1)}= B z,~~B_1: z\mapsto j = B_1z       \label{eq:j}
 \ee
 are generic linear
 mappings from the ${\bf\rm PGL(2,C)}$ group and $I_s$ is
 an  involution, the standard Cremona mapping with three simple $F$-points $O_\alpha\in
\{(1,0,0), (0,1,0), (0,0,1)\}$  and three principal lines
 $J_\alpha= (z:j_\alpha(z)=0)\in\{(z_1=0),(z_2=0),(z_3=0)\}$:
\bea
 I_s & : & z\mapsto z' \quad z'_1:z'_2:z'_3 = z_2z_3:z_1z_3:z_1z_2.
 \label{eq:stand}\\
 I_s & : & J_\beta\to O^{(-1)}_\beta,~
 O^{(-1)}_\beta\in((1,0,0),(0,1,0),(0,0,1)),\label{eq:stand1}\\
 I^{(-1)}_s & : & z'\to z \quad z_1:z_2:z_3=z_2'z_3':z_1'z_3':z_1'z_2',\label{eq:stand-1}\\
 I^{(-1)}_s &: & J^{(-1)}_\alpha\to O_\alpha,~
 O_\alpha\in\{(1,0,0),(0,1,0),(0,0,1)\}. \label{eq:stand-11}
 \eea
  In the
triangular frame of reference \rf{eq:j} mapping \rf{eq:quad} takes a very simple form
\bea
 \Phi_2 :
j(z)\mapsto j^{(-1)}(z') \quad & & j^{(-1)}_1(z'):j^{(-1)}_2(z'):j^{(-1)}_3(z')=\nonumber \\
& & j_2(z)j_3(z):j_1(z)j_3(z):j_1(z)j_2(z). \label{eq:quadj}
 \eea
  The mapping $\Phi_2$ \rf{eq:quad} is specialized if two or three $F$-points are adjacent
 or infinitely near \citep{iskovskikh.reid-91} and has, respectively, the following forms:
 \bea
  \Phi_{2a} & \equiv & B^{-1}\comp I_a\comp B_1,~~I_a:
z\mapsto z'~~z'_1:z'_2:z'_3 =
 z_2^2:z_1z_2:z_1z_3,   \label{eq:quada}\\
\Phi_{2b} & \equiv & B^{-1}\comp I_b\comp B_1,~~I_b: z\mapsto z'~~z'_1:z'_2:z'_3 =
 z_1^2:z_1z_2:(z_2^2-z_1z_3),   \label{eq:quadb}
\eea
 moreover, involutions $I_a, I_b$ from \rf{eq:quada},
\rf{eq:quadb} can
 be resolved as a composition of two or four, but not fewer, general mappings
 \rf{eq:quad}, respectively (see \citep{hudson-27}, chapter III, pp. 35,37).
    Any two members of the net \rf{eq:quada} touch
one another and have a fixed common tangent $j_1 \equiv z_1=0$, but ones of
 the net \rf{eq:quadb} have a fixed common tangent $j\equiv z_1 $ and osculate
 a fixed conic $ z_2^2-z_1z_3$. These tangency conditions are simulated
by two or three infinitely near points, so as  equations
 \rf{eq:f6}-\rf{eq:f11}
remain correct.}

\end{defin}

\subsection{\bf Properties of Birational Mapping}\label{properties}

\begin{defin} {\bf Properties.}\label{prop}
{\rm The one-to-one correspondence  for direct $\Phi_n$ and inverse $
\Phi^{-1}_n$ mappings does not hold only at indeterminacy or
fundamental points ($F$-points)
  $O_\alpha \in {\bf O},$ $O^{(-1)}_\beta \in {\bf O^{(-1)}},$
  $\alpha, \beta=(1,2,\ldots,\sigma),$
i.e., common basis points of multiplicities $i_\alpha, i^{(-1)}_\beta$ for functions $\phi_k(z),
\phi^{(-1)}_k(z), k=(1,2,3), $ and the associated
 linear systems $\phi$ \rf{eq:fa7} and $\phi^{(-1)}$ \rf{eq:fa8}, respectively,
 and on principal or exceptional curves or exceptional divisors
 $J_\beta, J^{(-1)}_\alpha,  ~ \alpha, \beta=(1,2,\ldots,\sigma)$,
\be
 J_\beta\defeq\{z:j_\beta(z)=0\},~~J^{(-1)}_\alpha\defeq\{z:j^{(-1)}_\alpha(z)=0\},~~
\alpha,\beta=(1, \cdots,\sigma), \label{eq:fa9}
 \ee
  where  $j_\beta, j^{(-1)}_\alpha $ are homogeneous polynomials in $z$ of degrees
 $i^{(-1)}_\beta, i_\alpha $, respectively, moreover, the points $O_\alpha, O^{(-1)}_\beta$
blow up into the curves $J^{(-1)}_\alpha, J_\beta$ of degrees $i_\alpha, i^{(-1)}_\beta$ and the curves
$J^{(-1)}_\alpha, J_\beta$ blow down into the points $O_\alpha, O^{(-1)}_\beta$,
 \bea
O_\alpha & \stackrel{\rightarrow}{\leftarrow} & J^{(-1)}_\alpha,\quad
 {\rm deg} J^{(-1)}_\alpha=i_\alpha,\label{eq:OJ'}\\
O^{(-1)}_\beta & \stackrel{\rightarrow}{\leftarrow} & J_\beta,~~~~\quad{\rm deg} J_\beta= i^{(-1)}_\beta,
\label{eq:O'J}
\eea
 respectively (see the
concept of $\sigma$-process of blowing up of singularities in the theory of ordinary differential
equations \citep{arnold-88} and the Kodaira theorem in the algebraic geometry
\citep{griffiths.harris-78}). $\triangleleft$}
\end{defin}

\begin{defin} \label{inf.near.points}
 {\rm
  A fundamental point is called ordinary  if
at this point there are no any additional tangency conditions.
 In special cases of non-ordinary (extra-ordinary) F-points
 (see Definition \ref{non.ord.point}),
tangency conditions of any two
members of the associated linear systems are expressed as multiplicities of
infinitely near points \citep{iskovskikh.reid-91}, or adjoint points in
 the terminology of \citep{hudson-27}.  Each $r$-fold non-ordinary point can
  be represented by the system of infinitely near ordinary points  and be transformed into ordinary
  ones using the technics of resolution of singularities of plane curves (see \citep{hudson-27},
   Chapter VII, Theorem \ref{equiv.map} below and examples :  \ref{ex1}, \ref{ex2} in
 Appendices A, B).
}
\end{defin}

\begin{theor}{\bf Jacobian.} \label{Jacobian}
{\rm
The Jacobian $J$ of the mapping $\Phi_n$ equals
\be
J = \left\|\frac {\partial \phi_k}{\partial z_i}\right\| \sim
    \prod_{\alpha = 1}^{\sigma}j_\alpha.\qquad {\rm deg} J= 3 n-3
    \label{eq:Jac}
\ee
 The formula corresponds to a birational map with ordinary $F-$
points but in the  case of non-ordinary \\( infinitely near) points it remains correct if we assign
to the infinitely near $F-$ points the same factors $j_\alpha$ with multiplicities in accordance
with a characteristic of the map (see below Definition \ref{characteristic} and
also  examples \ref{ex1}, \ref{ex2} in Appendix A).
 The determination of the Jacobian is a very simple way to find
the principal curves.  The principal curves $J_\alpha (J^{(-1)}_\beta)$ intersect each other only
at fundamental points $O_\alpha (O^{(-1)}_\beta)$. }
\end{theor}

\begin{rem} {\bf Characteristic.} \label{characteristic}
{\rm
 The set of numbers
 $$char(\Phi_n)=\{n;i_1,i_2,\ldots,i_\sigma\},\quad,
 i_1\geq i_2\geq\cdots\geq i_\sigma,
 $$
  where $i_\alpha $ are  the multiplicities of all indeterminacy points $O_\alpha$ of the
   mapping $\Phi_n$,
   including infinitely near ones, is called the characteristic of
 mapping $\Phi_n$. We shall denote the infinitely near $F-$points by
 the star: $i^{*}_\alpha$ and $O^{*}_\alpha$.

Next  in simplicity after quadratic birational map with
 $char(\Phi_2)=\{2; 1,1,1\}$ is a cubic
 map with $char(\Phi_3)=\{3; 2, 1, 1, 1, 1\}$ and then two quartic maps with
   $char(\Phi_4)=\{4; 2, 2, 2, 1, 1, 1\}$ and $\{4; 3, 1, 1, 1, 1,\\ 1, 1\}$.

   The general
 mapping with a given characteristic depends on $2\sigma +8$ parameters.

$\triangleleft$
}
\end{rem}


\begin{rem}{\bf Characteristic numbers.} \label{numbers}
{\rm
Let $i^{(-1)}_{\beta\alpha}$ be the multiplicity of curve
  $J^{(-1)}_\alpha$ at point $O^{(-1)}_\beta$
and $i_{\alpha\beta}$ be that of curve $J_\beta$ at $O_\alpha$. Then  we have
the equality $i_{\alpha\beta}=i^{(-1)}_{\beta\alpha}$ and the following relations
between
numbers $i_\alpha, i^{(-1)}_\beta,i_{\alpha\beta}$, expressing certain geometrical
 facts (summing in the left column over $\alpha $ and in the right one over $\beta $ from 1 to
$\sigma$):
\bea
\sum i_\alpha~~~~~ &=&~~~~~~~~~~~~~~ 3(n-1), ~~\sum i^{(-1)}_\beta~~~~~~~~~ = ~~~~~~ 3(n-1),
  \label{eq:f6}  \\
\sum i_{\alpha}^2~~~~~ &=&~~~~~~~~~~~~~~~~ n^2-1, ~~~\sum {i^{(-1)}_{\beta}}^2~~~~~~~ = ~~~~~~~~~
 n^2-1,   \label{eq:f7}  \\
\sum i_{\alpha\beta}~~~~ &=&~~~~~~~~~~~~ 3i^{(-1)}_\beta-1, ~~~
\sum i_{\alpha\beta}~~~~~~~~~~ = ~~~~~~~~ 3i_
\alpha-1,   \label{eq:f8} \\
\sum i_\alpha i_{\alpha\beta}~ &=&~~~~~~~~~~~~~~~~~ i^{(-1)}_\beta
n, ~~~ \sum i^{(-1)}_\beta
 i_{\alpha\beta}~~~~ = ~~~~~~~~~~~~ i_\alpha n,   \label{eq:f9} \\
\sum i_{\alpha\beta} i_{\alpha\gamma} &=&~~~~~ i^{(-1)}_\beta
i^{(-1)}_\gamma + \delta_{\beta\gamma},~~~\sum
i_{\alpha\beta}i_{\gamma\beta}~~~~~~ = ~~~~ i_\alpha i_\gamma +
\delta_{\alpha\gamma}.    \label{eq:f11}
 \eea
 The conditions
\rf{eq:f6}, \rf{eq:f7}  mean that the associated linear systems
\rf{eq:fa7}, \rf{eq:fa8} have the grade  $ D=1$, the genus $p=0$,
 the dimension $k=2$, and the superabundance $s=0$. The conditions \rf{eq:f8},
\rf{eq:f11} provide  rationality of the curves $J_\beta, J^{(-1)}_\alpha$
 \rf{eq:fa9}, and
that their degrees are $i^{(-1)}_\beta$ and $i_\alpha$, respectively.
In the case of non-ordinary $F$-points the total number of distinct $F$-points
   need not be the same for the direct \rf{eq:fin2} and inverse \rf{eq:fin2'}
mappings. In the special cases, if at some $\beta$   in the left parts of
 equations \rf{eq:f8} and \rf{eq:f11} (at $\beta=\gamma  $) $j_\beta$ breaks up into $1\leq\nu
\leq i^{(-1)}_\beta$ components,  then the left parts of these  equations must be
 replaced by
\be
\sum i_{\alpha\beta} = 3i^{(-1)}_\beta-\nu,~~~~\sum i^2_{\alpha\beta} =
 i^{{(-1)}^2}_\beta +\nu. \label{eq:brekjbeta}
\ee
 Analogous changes in the right parts of equations \rf{eq:f8} and
 \rf{eq:f11} at some $\alpha$ must
be made.

The upper limit for $\sigma$  of the total number of $F$-points is given by
the following  formula:
$$
\sigma \leq 2n-1, \qquad \mbox{if}~~ n>1.
$$
 $\triangleleft$
 }
\end{rem}

\subsection{\bf Behaviour of algebraic curves and LSAC
under the action of the birational map \rf{eq:fin2}} \label{curve.lsc}

\begin{rem}  \label{curve}
{\rm
 Consider  properties of a general curve $f_\mu(z')=0 $ of degree $\mu$
 under the mapping
 \rf{eq:fin2}. By map \rf{eq:fin2}, the curve $f_\mu(z')$ is mapped into  the curve
 $f_\mu(\phi(z))=f'_{\mu'}(z)$ of degree  $\mu'=\mu n$; moreover, every point
 $O_\alpha$ which is $i_\alpha$-fold on $\phi(z)$ is $\mu i_\alpha$-fold on
 $f'_{\mu'}$.
 If $f_\mu(z')$ has multiplicities $\gamma^{(-1)}_\beta$ at points $O^{(-1)}_\beta$,
 then (${\rm deg}(j_\beta) \equiv i^{(-1)}_\beta $)
\be
  f_\mu(z')=f'_{\mu'}(z)\prod_{\beta=1}^{\sigma} j_{\beta}^{\gamma^{(-1)}_\beta},
  ~~~\mu'=\mu n -\sum_{\beta=1}^{\sigma}\gamma^{(-1)}_\beta i^{(-1)}_\beta;
                                      \label{eq:f12}
\ee
 moreover, $f'_{\mu'}$ has multiplicities $\gamma'_\alpha$ at $O_\alpha$
 (see the meaning of $i_{\alpha \beta}$ in Remark \ref{numbers}):
\be
  \gamma'_\alpha =\mu i_\alpha- \sum_{\beta=1}^{\sigma}
                 i_{\alpha\beta}\gamma^{(-1)}_\beta.   \label{eq:f13}
\ee

If $f_\mu(z)=0$ is a general curve of a linear system of curves of dimension $k$  but $f'_{\mu'}$
is its image under the map $\Phi_n$ \rf{eq:fin2}, that genus $p$ and dimension $k$ of the LSAC
 are invariants.
}
\end{rem}

\subsection{\bf Birational equivalence and  resolution of singularities of plane curves}
 \label{bir.equiv.res.sing}

Let us introduce the definition of birationally equivalent mapping in ${\bf\rm CP}^2$.
\begin{defin} \label{bir.equiv1}
{\rm
A mapping $\Phi_n \in {\bf\rm Bir CP}^2, z\to z'\sim\phi(z),$
is birationally equivalent (or conjugated in terminology of Hudson's book)
 to a mapping $U: y\to y', y'\sim u(y)$  $y, y' \in {\bf\rm CP}^2,$
if there exists a
birational mapping $V_m:z\to y, y\sim v(z) $ of degree $m$
 such that $\Phi_n \equiv V_m^{-1}\comp U \comp V_m$ .
 }
\end{defin}

Due to the standard method of resolution of singularities of plane curves
(see \citep{hudson-27}, chapter VII, p.129)  the following problems can be
 solved
by applying a composition of the corresponding Cremona quadratic mappings:

\begin{enumerate}
\item
to transform any non-ordinary multiple point into a net of simple points;
\item
to resolve any non-ordinary multiple point into an equivalent set of ordinary
multiple points;
\item
to transform any algebraic curve into one having ordinary multiple points
only;
\item
to transform any linear system of algebraic curves into one having ordinary
base points only.
\end{enumerate}

As the consequence of the standard method of resolution of singularities of plane curves the following
 theorem is represented to be valid.

\begin{theor} \label{equiv.map}
{\rm Any mapping $\Phi_n$ \rf{eq:fin2} with non-ordinary $F$-points (see
Definition \ref{non.ord.point}) by the
 corresponding transformation of birational equivalence (see Definition  \ref{bir.equiv1})
 is transformed
into  some mapping $\Phi_{n'}=V^{-1}\comp \Phi_n\comp V$ with only ordinary
$F$-points where degree of this mapping $n'\geq n$ but the mapping $V$ is a composition of
 a necessary (for resolution of all infinitely near points) number of quadratic mappings.
   (see Definition \ref{quadratic})
  ( see  example in Appendix B) ~~$\triangleleft$
}
\end{theor}

For illustration of this Section we set the examples of quadratic and cubic maps with ordinary
 and non-ordinary indeterminacy points in Appendix A but the example of using Theorem
 \ref{equiv.map} is in Appendix B.

Below we shall deal with mappings having only ordinary indeterminacy points supposing that
maps with non-ordinary ones was previously replaced by birationally equivalent maps with
the help of Theorem \ref{equiv.map}.


\section{Decomposition of the set of indeterminacy points} \label{decomp}
Let us consider orbits of indeterminacy points $O_{\alpha}\in {\bf\rm O} \quad \mbox{and}\quad
O^{(-1)}_\beta\in {\bf\rm O}^{(-1)}$ relative to the action of the inverse $\Phi^{-1}_n$
 \rf{eq:fin2'} and the direct map $\Phi_n$ \rf{eq:fin2}, respectively, and let us introduce the following definitions.

\begin{defin}  \label{orbit}
{\rm The orbit ${\bf\cal O}_z$ of a point $z$ with respect
 to  the mapping $\Phi^{-1}_n$ \rf{eq:fin2'} is the set of
 points ${\cal O}^k_z = \Phi^{-k}_n(z) = (\Phi^{-1}_n)^{k}(z),~
 k\in {\bf\rm Z}^{+}$,
 where ${\bf\rm Z}^{+}$ is the set of non-negative integers.
The orbit ${\bf\cal O}^{(-1)}_z$ of a point z with respect to $\Phi_n$
\rf{eq:fin2} is defined analogously, $ {\cal O}^{(-1)k}_z = \Phi^{k}_n(z) =
(\Phi_n)^{k}(z)$ and $
 \Phi^k_n(z)\defeq \Phi_n(\Phi_n(\ldots(z)\ldots)), \Phi^{-k}_n(z)\defeq
 \Phi^{-1}_n(\Phi^{-1}_n(\ldots(z)\ldots))
$ (see, for example, \citep{DS-1-88a}, \citep{arnold-88}).}
\end{defin}

\begin{defin} \label{cycle}
{\rm If the number $k$ of the points of the orbit ${\bf\cal O}_z$
  of the point $z$ with respect to  the mapping $\Phi^{-1}_n$ \rf{eq:fin2'}
is finite, where non-negative integer $k$ is a minimal integer defined by
 the condition
$$
{\cal O}^k_z = \Phi^{-k}_n(z) = z,~~ k\in {\bf\rm Z}^{+},
$$
then the periodic points $ \bigl(\Phi^{-m}_n(z)\bigr),$
$m=\bigl(0,1,\cdots k-1\bigr)$ form the set ${\bf\cal O}^{(cycle)}_z$-- a cycle of
index $k$ or period $k$ of the mapping $\Phi^{-1}_n(z)$,
 but $z$ is a fixed point of the mapping
$\Phi^{-k}_n(z)$. A cycle of index $k$ of the mapping $\Phi_n(z)$
\rf{eq:fin2} is defined similarly with  help of the changes: ${\bf\cal O}_z \mapsto
{\bf\cal O}^{(-1)}_z$ and $ \Phi^{-m}_n(z) \mapsto \Phi^{m}_n(z),$ (see
\citep{DS-1-88a}).
}
\end{defin}

Let us introduce the  notion of a tail of the cycle.

\begin{defin} \label{tail}
{\rm
Let us call a subset ${\bf\cal O}^{(tail)}_z$ of the set ${\bf\cal O}_z$ a  tail of the length $l$ of
the cycle ${\bf\cal O}^{(cycle)}_y$
 where non-negative integer $l$ is a minimal integer
 defined by the condition
$$
{\bf\cal O}^{(tail)}_z\mbox{ of length}~~ l: \{y = \Phi^{-l}_n(z), \mbox{
where}~ y \in{\bf {\cal O}}^{(cycle)}_y\},~~l
 \in {\bf\rm Z}^{+},
$$
so that the point $z$ is the beginning of the tail, but the point $y$ is  the beginning of
the cycle and does not belong to the tail.
}
\end{defin}

\begin{defin}{\bf Decomposition}  \label{decomposition}

{\rm Let $\Phi_n$ \rf{eq:fin2} be a mapping of characteristic
$n;i_1,i_2,\cdots, i_\sigma$  and $\Phi^{-1}_n$ \rf{eq:fin2'} be the
inverse mapping (see Definition-Theorem \ref{map}, Remark
\ref{characteristic}). Define \citep{rerikh-98a}
 the decomposition of the sets ${\bf O, O^{(-1)}}$
 of fundamental points $O_\alpha, O^{(-1)}_\beta$ of these mappings
as follows:
\bea
{\bf O} & \equiv & {\bf O}^{(rest)}\cup{\bf O}^{(int)}, \quad {\bf O^{(-1)}} \equiv{\bf
O^{(-1)}}^{(rest)}\cup{\bf O^{(-1)}}^{(int)}, \label{eq:f19} \\
{\bf O}^{(rest)} & \equiv & {\bf O}^{(cycle)}\cup{\bf O}^{(tails)}\cup{\bf
O}^{(inf)}, \nonumber \\
{\bf O^{(-1)}}^{(rest)} & \equiv &
 {\bf O^{(-1)}}^{(cycle)}\cup{\bf O^{(-1)}}^{(tails)}\cup{\bf O^{(-1)}}^{(inf)}. \label{eq:f19a}
 \eea
  Here: ${\bf O}^{(inf)}$ is a subset of
fundamental points $O_\alpha$ with infinite orbits
\be
 {\bf\rm
O}^{(inf)}~~~:[\Phi^{-k}_n({\bf O})\cap[{\bf O}^{(-1)}\cup
{\bf O}]]=\emptyset \quad\mbox{at}~\forall k ~ 1\leq k <\infty;
\label{eq:oinf}
\ee
 ${\bf O}^{(cycle)}$ is a subset of fundamental points $O_\alpha$ having
cyclic orbits ${\cal O}^{m}_z, z \in {\bf O}$,  of index $m_\alpha$;
 ${\bf O}^{(tails)}$ is a subset of fundamental points $O_\alpha$ belonging to the tails
 of the orbits of the subset ${\bf O}^{(cycle)}$, to the tails
  of the orbits of the subset ${\bf O}^{(int)}$ and to the tails of the orbits
  of the subset ${\bf O}^{(inf)}$.

  The subsets ${\bf O}^{(-1)(int)}$ and ${\bf O}^{(int)}$
of the sets ${\bf O^{(-1)}}$ and ${\bf O}$ are defined below.
  }
  \begin{defin} \label{oint}
{\rm The subsets ${\bf O^{(-1)}}^{(int)}$ and ${\bf O}^{(int)}$
of the sets ${\bf O^{(-1)}}$ and ${\bf O}$ are under
construction in the following manner.~ Let  ${\Phi^{-k}_n}({\bf O}(k)),
~ {\Phi^{k}_n}({\bf O^{(-1)}}(k))$,
$ k \geq 0,~ \Phi^{-k}_n|_{k=0}\equiv {\rm id}$
be  $k^{th}-$iterations of
 punctual sets ${\bf O}(k),~{\bf O^{(-1)}}(k),~ {\bf O}(0)$ $\equiv
 {\bf O},\quad {\bf O^{(-1)}}(0)\equiv {\bf O^{(-1)}}$  under the action of
   inverse and direct maps
$\Phi^{-1}_n, \Phi_n$. Let us introduce the punctual set
\bea
 {\bf O}^{(int)}(k) & \stackrel{\rm def}{=} &  \Phi^{-k}_n({\bf O}(k))\bigcap {\bf O}^{(-1)}(k),\nonumber\\
{\bf O}^{(-1)(int)}(k) & \stackrel{\rm def}{=} &  \Phi^{k}_n({\bf O^{(-1)}}(k))\bigcap {\bf O}(k),
  \quad k\geq 0 \label{eq:oint(k)}
\eea
 and define  the construction of the set ${\bf O}(k)\quad\mbox{and}\quad {\bf O}^{(-1)}(k)\quad
 \mbox{at}\quad k\geq 1 $ in
the following manner:
 \bea
  {\bf O}(k) & \stackrel{\rm def}{=} &
{\bf O}(k-1)/{\bf O}^{(int)}(k-1), \nonumber \\
{\bf O}^{(-1)}(k) & \stackrel{\rm def}{=} & {\bf O}^{(-1)}(k-1)/{\bf O}^{(-1)(int)}(k-1), ~~ k\geq 1. \label{eq:O(k)}
 \eea
  Then  subsets
${\bf O}^{(int)}$ and ${\bf O}^{(-1)(int)}$ are defined  with the help of \rf{eq:oint(k)} and
\rf{eq:O(k)} by
\be
 {\bf O}^{(int)}\stackrel{\rm
def}{=}\bigcup_{k=0}^{k=m}{\bf O}^{(int)}(k),\quad
{\bf O}^{(-1)(int)}\stackrel{\rm
def}{=}\bigcup_{k=0}^{k=m}{\bf O}^{(-1)(int)}(k),
 \label{eq:Oint}
\ee
 where positive integer $m$ is defined by the condition
 $$
{\bf O}(m+1)/{\bf O}^{(rest)}\equiv \emptyset,\quad
{\bf O}^{(-1)}(m+1)/{\bf O}^{(-1)(rest)}\equiv \emptyset,
 $$
 where there subsets ${\bf O}^{(rest)}, {\bf O}^{(-1)(rest)}$ are defined above by
\rf{eq:f19a}. Let be
 \be
  \# {\bf O}^{(int)}=\sigma_1.
\label{eq:sigma1}
\ee
 The above construction
 establishes the one-to-one correspondence between $\sigma_1$
pairs of
 equivalent
indeterminacy points of the subsets ${\bf\rm O}^{(int)}$ and
${\bf\rm O}^{(-1) (int)}$ constructed above as follows:
\bea
& & \Phi^{-m_j}(O_{\alpha_j}) \equiv  O^{(-1)}_{\beta_j}
,\quad j=(1,\cdots,\sigma_1), \label{eq:equiv.pairs1} \\
& & O_{\alpha_j} \in {\bf\rm O}^{(int)},\quad O^{(-1)}_{\beta_j}\in {\bf\rm
O}^{(-1)(int)},\quad \alpha_j ~~\mbox{and}~~ \beta_j  \in  (1,\cdots,\sigma),
\label{eq:equiv.pairs2}\\
& & m_j  \in  (m_1,\cdots,m_{\sigma_1})\quad 0\leq m_1\leq\cdots,\leq
m_{\sigma_1}, \label{eq:set.m}
 \eea
 where nonnegative integers $m_j$ are lengths of the orbits of points $O_{\alpha_j}$ but integer
 $m$ in equation \rf{eq:Oint} is equal to $m_{\sigma_1}$.
} $\triangleleft$
\begin{rem} \label{integrable.surfaces. in.parameterical.space}
{
\rm~~~ Since the coordinates of the indeterminacy points $O_{\alpha_j}, O^{(-1)}_{\beta_j}$,
 $j=(1,2,\cdots,\sigma_1)$
are functions of $2\sigma+8$ parameters, equations \rf{eq:equiv.pairs1}-\rf{eq:set.m}
define  in the space of $(2\sigma+8)$ parameters $2\sigma_1$ subvarieties of  dimension
$2\sigma+8-2\sigma_1$. If the A.P.~Veselov conjecture is true, integrable mappings correspond
to these subvarieties.
}
\end{rem}
$\triangleleft$
\end{defin}
 \end{defin}

\section{Dynamics of a generic birational mapping.\\
Difference equation for the Arnold complexity.} \label{dynamics.diff.eq}

 Theorem \ref{dynamics}  defines the dynamics
 of birational mapping and the difference equation for $d(k)$.

\begin{theor} \label{dynamics}

~~~{\rm Let $d(k)$ be  the degree of the mapping $\Phi_n^k$
\be
\Phi^k_n\quad :\quad z\to z'\qquad z'_1:z'_2:z'_3 = \phi^{(k)}_1(z):\phi^{(k)}_2(z):
\phi^{(k)}_3(z), \label{eq:f21}
\ee
 the $k^{th}$ iteration of the mapping
 $\Phi_n$ \rf{eq:fin2} of characteristic $char(\Phi_n)=\{n,i_1,\cdots,\\ i_{\sigma}\}$
 , $O_{\alpha_j}$ and $O^{(-1)}_{\beta_j}, j=(1,\cdots,\sigma_1)$ be $\sigma_1$
 pairs of equivalent indeterminacy points of the  subsets ${\bf\rm O}^{(int)}\quad
 \mbox{and}~ {\bf\rm O}^{(-1)(int)}$ defined by relations \rf{eq:equiv.pairs1}, \rf{eq:equiv.pairs2}
 and \rf{eq:set.m} (see Section \ref{decomp}, Definitions \ref{decomposition} and \ref{oint}).
 Let also  $\gamma_{\alpha_j}(k)$
 be common multiplicities of the curves $\{\phi_i^{(k)}(z)=0, i=(1,2,3)\}$
 and the general curve  of the linear system  $\phi^{(k)}_\mu(z)=\{\sum_{i=1}^{i=3}c_i \phi_i^{(k)}(z)=0, \quad
 \forall c_i \in{\bf\rm C}, \}$ of degree $\mu=d(k)$
  at indeterminacy points $O_{\alpha_j}$ of the
 direct mapping \rf{eq:fin2}
  (we assume that all F-points are already ordinary after the birational equivalence
  transformation --see  Theorem \ref{equiv.map}). Then the dynamics of
 the mapping $\Phi_n$ \rf{eq:fin2} (see  Definition
\ref{map}, Remarks \ref{characteristic} and \ref{numbers})
 is completely determined by the following set of difference equations:
\bea
  d(k) & = & nd(k-1)-\sum^{\sigma_1}_{l=1} i^{(-1)}_{\beta_l}\gamma_{\alpha_l}(k-m_l-1),
     \label{eq:f22}\\
  \gamma_{\alpha_j}(k) & = & i_{\alpha_j}d(k-1) -
                         \sum^{\sigma 1}_{l=1} i_{\alpha_j\beta_l}\gamma_
             {\alpha_l}(k-m_l-1), ~ j=1,\cdots,\sigma_1,
                                                  \label{eq:f23}\\
\gamma_{\alpha}(k) & = & i_{\alpha}d(k-1) -
                         \sum^{\sigma 1}_{l=1} i_{\alpha\beta_l}\gamma_
             {\alpha_l}(k-m_l-1), \qquad \alpha\neq\alpha_j,
                                                  \label{eq:f23a}
\eea
moreover,
\be
d(0) = 1,~~d(1)=n,~~\gamma_\alpha(1)=i_\alpha,~~\gamma_\alpha(k)= 0~~
 \mbox{for}~~ k\leq 0.    \label{eq:f24}
 \ee

The secular equation corresponding to the set of difference  equations \rf{eq:f22}-\rf{eq:f23} is
\be
\det(\Lambda)=\lambda^m + \sum_{i=0}^{m-1}a_i\lambda^i = 0, \label{eq:sec.eq}
\ee
where integer m is
 \be
  m =m_1+m_2+\cdots+m_{\sigma_1}+\sigma_1+1,  \label{eq:m}
 \ee
integers $a_i$ are coefficients of expansion in power series of $\det(\Lambda)$ in $\lambda$
 and the matrix $\Lambda$ is
\be
\Lambda=\left(
\begin{array}{ccccc}
  \lambda-n, & i^{(-1)}_{\beta_1},  & \cdots, & i^{(-1)}_{\beta_{\sigma_1}}\\
  -i_{\alpha_1}, & (\lambda^{m_1+1}+i_{\alpha_1\beta_1}),  & \cdots, &
  i_{\alpha_1\beta_{\sigma_1}} \label{eq:lambda}\\
    \vdots & \vdots & \vdots \vdots \vdots & \vdots \\
   -i_{\alpha_{\sigma_1}}, & i_{\alpha_{\sigma_1}\beta_1}, &
   \cdots, &
   (\lambda^{m_{\sigma_1}+1}+i_{\alpha_{\sigma_1}\beta_{\sigma_1}})
\end{array}
\right).
\ee
 The linear difference equation for $d(k)$ corresponding to the secular equation \rf{eq:sec.eq}  has
the form
\be
 d(k+m) + \sum_{i=0}^{m-1}a_id(k+i)=0,
\label{eq:f39}
\ee
where $a_i$ are the same integers as in equation \rf{eq:sec.eq}.

According to  a general theory of linear difference equations with constant coefficients
 \citep{gel'fond-67} (Chapter V), the solution of equation \rf{eq:f39} has
 the form
 \be
d(k)=\sum_{i=1}^{l}\lambda_i^k(\sum_{j=0}^{s_i-1}c_{ij}k^j), \label{eq:f41}
\ee
where $\lambda_1,...,\lambda_l$ are multiple roots of equation \rf{eq:sec.eq}
with multiplicities $s_1,...,s_l,$ \- $s_1+s_2+\cdots+s_l=m, $ and $c_{ij}$
are arbitrary constants to be determined from $m$ initial values
$d(1)=n, d(2), ..., d(m)$ obtained with the help of equations \rf{eq:f22}-
\rf{eq:f24}.

It is obvious that the condition
\be
|\lambda_i|=1 \quad \forall i\in(1,2,\cdots,l) \label{eq:pol.bound}
\ee
is sufficient for the polynomial boundedness of the growth of $d(k)$ with $k$.

\begin{rem}\label{equal.some.m_i}
{\rm
If some integers $m_i$ are equal and, moreover, are fulfilled corresponding conditions for the coefficients
 at the terms $\gamma_{\alpha_i}(k)$ in equations \rf{eq:f22}-\rf{eq:f23},
 we can decrease  the order of the system of difference equations
 \rf{eq:f22}-\rf{eq:f23} and, as result, the order of the matrix $\Lambda$.
 }
 \end{rem}
$\triangleleft $ }
\end{theor}
Proof.~~~ Let us prove the theorem by an induction method. Let us consider
 the map $\Phi^k_n$  \rf{eq:f21}-the
$k^{th} $ iteration of the mapping $\Phi_n$ \rf{eq:fin2} as an iteration of the map
$\Phi^{k-1}_n$  and let us  consider  the
transformation of a general curve of a linear system of the curves
$\phi^{(k-1)}_\mu(z)=\sum_{i=1}^{i=3} c_i\phi^{(k-1)}_i(z)=0$ of degree $\mu=d(k-1)$
 by the action of the
mapping $\Phi_n$ \rf{eq:fin2}.  Let $\gamma^{(-1)}_{\beta}(k-1)$
 be common multiplicities of the curves  $\{\phi^{(k-1)}_i(z)=0, \quad \forall i\in (1,2,3)\}$ and of the general curve
  of the linear system $\{\phi^{(k-1)}_\mu(z)=0\}$  at the indeterminacy  points $O^{(-1)}_{\beta}$
  of the inverse map $\Phi^{(-1)}_n$ \rf{eq:fin2'}. Then, according to Subsection \ref{curve.lsc},
  Remark \ref{curve}, \rf{eq:f12} and \rf{eq:f13},
   we have
  \bea
  \phi^{(k-1)}_\mu(\phi(z)) & = & \phi^{(k)}_{\mu'}(z)\prod_{\beta=1}^{\sigma} j_{\beta}^
  {\gamma^{(-1)}_\beta(k-1)}(z), \label{eq:f42}  \\
  \mu' & = & \mu n -\sum_{\beta=1}^{\sigma}i^{(-1)}_\beta \gamma^{(-1)}_\beta(k-1) ,\label{eq:f43} \\
  \gamma_\alpha(k) & = & \mu i_\alpha- \sum_{\beta=1}^{\sigma}i_{\alpha\beta}
  \gamma^{(-1)}_\beta(k-1),   \label{eq:f44}
  \eea
  where $\phi^{(k)}_{\mu'}(z) = \sum_{i=1}^{i=3} c_i\phi^{(k)}_i(z)$ is a general curve of a linear
   system of curves of degree $\mu'=d(k)$ associated with the map $\Phi^k_n$ but $\gamma^{(k)}_\alpha$
    are its multiplicities at the points $O_\alpha$.
   Since the linear system of the curves $\{\phi^{(k-1)}_\mu(z)=0\}$ is completely defined by
   its basis set, the difference of  values $\gamma^{(-1)}_\beta(k-1)$ from zero means that the set
   ${\bf\rm O}^{(k-1)}\bigcap {\bf\rm O}^{(-1)} \neq \emptyset$ where ${\bf\rm O}^{(k-1)}$ is
   the set of indeterminacy points the mapping $\Phi^{k-1}_n$. It is obvious that the set
   ${\bf\rm O}^{(k-1)}$ is equal to (${\bf\rm O}^{(1)}\equiv {\bf\rm O}$)
   \be
    {\bf\rm O}^{(k-1)}=\bigcup_{l=0}^{l=k-2}\Phi^{-l}_n({\bf\rm O}(l)), \label{eq:f45}
    \ee
   where the set ${\bf\rm O}(l)$ is defined by equation \rf{eq:O(k)} (see Section \ref{decomp}).
   Let us decompose the set ${\bf\rm O}^{(k-1)}$ into two subsets
   \be
    {\bf\rm O}^{(k-1)}= {\bf\rm O}^{(k-1)(int)} \bigcup {\bf\rm O}^{(k-1)(rest)} \label{eq:f46}
    \ee
    related with the subsets ${\bf\rm O}^{(int)}$ and ${\bf\rm O}^{(rest)}$ in equation \rf{eq:f19}.
    Then the subset ${\bf\rm O}^{(k-1)(int)}$ is equal to
    \be
    {\bf\rm O}^{(k-1)(int)}=\bigcup_{j=1}^{j=\sigma_1}\bigcup_{l=0}^{l=\min(m_j,k-2)}
    \Phi^{-l}_n(O_{\alpha_j}), \label{eq:f47}
    \ee
    according to Section \ref{decomp}, Definitions \ref{decomposition}, \ref{oint} and
    \rf{eq:equiv.pairs1}, \rf{eq:equiv.pairs2} and \rf{eq:set.m}.

  Since, according to Section \ref{decomp},  the intersection of the subsets
  ${\bf\rm O}^{(k-1)(rest)}$ and ${\bf\rm O}^{(-1)}$ is empty, then $ {\bf\rm O}^{(k-1)}\bigcap {\bf\rm O}^{(-1)}
  = {\bf\rm O}^{(k-1)(int)}\bigcap{\bf\rm O}^{(-1)}$ and
  \be
    {\bf\rm O}^{(k-1)}\bigcap {\bf\rm O}^{(-1)}  =
\left\{\begin{array} {ccc}
               \bigcup_{j=1}^{j=\sigma_1}\Phi^{-m_j}_n(O_{\alpha_j})
                & =\bigcup_{j=1}^{j=\sigma_1}O^{(-1)}_{\beta_j},~& \forall (k-2)\geq m_j,\\
               \emptyset \quad
                &                                        & \forall (k-2)< m_j.
        \end{array}\right\} \label{eq:f48}
  \ee
Then, according to \rf{eq:f48}, we have for $\gamma^{(-1)}_\beta(k-1)=
{\rm mult}(\phi^{(k-1)}_\mu(z)|_{z=O^{(-1)}_\beta}$
 \be
 \gamma^{(-1)}_\beta(k-1)=
 \left\{ \begin{array}{cc}
   0,~ & \forall \beta\neq \beta_j,~ (\beta,\beta_j)\in(1,\cdots,\sigma) \\
   \gamma_{\alpha_j }(k-1-m_j),~ & \beta\equiv\beta_j,~ j\in(1,\cdots,\sigma_1)
   \end{array}\right\}, \label{eq:f49}
   \ee
   where
   \be
   \gamma_\alpha(k)=0,\quad k\leq 0,\qquad \gamma_\alpha(1)\defeq i_\alpha. \label{eq:f50}
   \ee
 At last, substituting \rf{eq:f49} into equations \rf{eq:f42}-\rf{eq:f44} and taking into account
  \rf{eq:f50} and $d(0)=1, d(1)=n$, we have equations \rf{eq:f22}-\rf{eq:f24}.
  Equations \rf{eq:f22}-\rf{eq:f24}  hold at $k=1$ and, therefore, at
 $\forall k>1$.
Let us prove the validity of equations \rf{eq:sec.eq}-\rf{eq:pol.bound}. There are two ways of
 obtaining them.
Let us transform  the system of $\sigma_1+1$ equations \rf{eq:f22} and \rf{eq:f23} at
 $\alpha_j=(\alpha_1, \alpha_2,...,\alpha_{\sigma_1})$ by changing
$k \rightarrow k+m_1 + 1$ to the following form $(j=1,\cdots,\sigma_1)$:
\bea
\sum_{j=1}^{\sigma_1}i^{(-1)}_{\beta_j}\gamma_{\alpha_j}(k+m_1-m_j) =
 nd(k+m_1)-d(k+m_1+1)
  \label{eq:f36} \\
\gamma_{\alpha_j}(k+m_1+1) + \sum_{l=1}^{\sigma_1}i_{\alpha_j \beta_j}\gamma_
{\alpha_l}(k+m_1-m_l)=d(k+m_1)i_{\alpha_j}.
 \label{eq:f37}
\eea

  The first way is to obtain a linear homogeneous difference equation
 \rf{eq:f39} for $d(k)$ excepting from the homogeneous system of $\sigma_1+1$ difference
 equations  \rf{eq:f36}, \rf{eq:f37} for unknown $d(k), \gamma_{\alpha_j}(k),$  $ j\in (1,2,
\cdots,\sigma_1)$ step by step
 $\gamma_{\alpha_1}(k), \gamma_{\alpha_2}(k)$, $\cdots$  and so on, until we do not obtain equation
 \rf{eq:f39} for $d(k)$.

 However, we can present a more direct method of obtaining equation \rf{eq:f39} for
the function $d(k)$ through finding the characteristic or secular equation immediately from the system
of equations \rf{eq:f36}-\rf{eq:f37} performing  substitution in them accordingly to
\be
d(k)=b_0 \lambda^k, \qquad \gamma_{\alpha_j}(k)=b_j \lambda^k, \quad j\in(1,2,\cdots,\sigma_1)
 \label{eq:subst}
\ee
where  $b_0, b_j, j\in(1,2, \cdots,\sigma_1)$ are unknown constants.

After the substitution \rf{eq:subst} the system of equations \rf{eq:f36}-\rf{eq:f37}
has the following matrix form:
\be
\Lambda D B=0,~~ D=diag(\lambda^{m_1},\lambda^{m_1-m_2}, \cdots, \lambda^{m_1-m_{\sigma_1}}),
~~ B=(b_0,b_1,\cdots,b_{\sigma_1}),\label{eq:f36',f37'}
\ee
where the matrix $\Lambda$ is defined by equation \rf{eq:lambda}.

The compatibility condition of the homogeneous system \rf{eq:f36',f37'} with respect to unknown
parameters $b_0, b_j, j\in(1,2,\cdots,\sigma_1)$ is a secular equation \rf{eq:sec.eq}
 where integers $a_i$, the same as in eq. \rf{eq:f39},  are
 the coefficients of expansion in power series of $\det (\Lambda)$ in $\lambda$ but
  integer $m$ is defined by \rf{eq:m}. Since the difference equation for $d(k)$ \rf{eq:f39} and
 the secular equation \rf{eq:sec.eq} are in one-to-one correspondence by the substitution
  $d(k)= \lambda^k$, we can uniquely reconstruct \rf{eq:f39} from \rf{eq:sec.eq}.

  At the end, according to the general theory of linear difference equations with
constant coefficients (see \citep{gel'fond-67}, chapter V), the
general solution of equation \rf{eq:f39}  is completely defined by
the spectrum of eigenvalues of the characteristic secular equation \rf{eq:sec.eq}
and this solution has the form \rf{eq:f41}. Remark that the method of obtaining a general
 solution of the system  of difference equations \rf{eq:f22}, \rf{eq:f23} offered above
 is fully analogous to usual practice of solving a system of linear differential equations of
 an order more one with constant coefficients ( see, for example, \citep{arnold-88c},
  the chapter 3,\S 25).
  $\triangleleft$

\section{Conclusion}

We believe in that all mappings with the Arnold complexity defined by equation \rf{eq:f41} and
the spectrum \rf{eq:pol.bound} of the secular equation \rf{eq:sec.eq} of the order $m$ \rf{eq:m}
with the matrix $\Lambda$ \rf{eq:lambda} are algebraically integrable ones and  are intended
to prove this theorem.

 Theorem \ref{dynamics} gives us
to possibility to generate, for example,  all integrable families of maps of degree $n=2$ in the parameter
 space of dimension $2\sigma+8-2\sigma_1=14-2\sigma_1$ being  stratified on  algebraic subvarieties
 in this space.
 We can present here some different interesting sets for $n=2, \sigma_1=3, \alpha_1=\beta_1,
\alpha_2=\beta_2, \alpha_3=\beta_3$:
 \bea
 m_1 & = & 0, m_2\neq m_3\neq 0, \quad \det(\Lambda)=(\lambda-1)^2(\lambda^{m_2+m_3+2}-1);
\label{eq:set1}\\
m_1 & = & 1, m_2=2, m_3=3, \quad \det(\Lambda)=(\lambda-1)^3(\lambda+1)\frac{(\lambda^9+1)}
{(\lambda^3+1)}; \label{eq:set2}\\
 m_1 & = & 1, m_2=2, m_3=4, \quad \det(\Lambda)=(\lambda-1)^3[\frac{(\lambda^{15}+1)
 (\lambda+1)}{(\lambda^5+1)(\lambda^3+1)}]; \label{eq:set3}\\
 m_1 & = & 1, m_2=2, m_3=5,~ \det(\Lambda)=(\lambda-1)^3(\lambda+1)(\lambda^5-1)
(\lambda^3-1); \label{eq:set4}\\
 m_1 & = & 0, m_2=m_3=m\geq 1, \quad \det(\Lambda)=(\lambda-1)^2(\lambda^{m+1}+1);
\label{eq:set5}\\
 m_1 & = & 1, m_2=m_3=2, \quad \det(\Lambda)=(\lambda-1)^2[\frac{\lambda^6+1}
 {\lambda^2+1}]; \label{eq:set6}\\
 m_1 & = & m_2=m_3=1,\quad \det(\Lambda)=(\lambda-1)\frac{\lambda^3+1}{\lambda+1};
 \label{eq:set7}\\
 m_1 & = & m_2=m_3=2, \det(\Lambda)=(\lambda-1)^3(\lambda+1). \label {eq:set8}
 \eea
 It isn't difficult to obtain polynomially abounded dependence $d(k)$ for sets \rf{eq:set1}--\rf{eq:set8}

 One presents itself interesting also to give classification of all integrable cubic maps.

\section{Acknowledgements}
The author is grateful to  V.A.~Iskovskikh, V.V.~Kozlov,
  A.N.~Parshin,  I.R.~Shafarevich and D.V.~Treschev
  for useful discussions and interest to the paper.

\section{ Appendix A. Examples of maps and FEs.} \label{examples}
\begin{enumerate}

\item{\bf FE of the paper.}\label{ex1}
\be
 y(w+2)=\frac{y(w+1)(\lambda y(w+1)+d y(w))}{y(w)}, \label{eq:ex01}
  \ee
 Supposing $y(w+2)=y'_2,~ y(w+1)=y'_1=y_2,~ y(w)=y_1$  we have  the map
  $y\mapsto y':{\rm\bf C}^2\to{\rm\bf C}^2$ and then changing
 $y\mapsto z: y_i=z_i/z_3$ we obtain the map $z\mapsto z':{\rm\bf CP}^2\to{\rm\bf CP}^2$:
 \bea
 \Phi_2 & : & ~ z'_1:z'_2:z'_3=z_2 z_1:z_2(\lambda z_2+d z_1):z_1
 z_3, \label{eq:ex1}\\
 \Phi^{(-1)}_2 & : & ~ z_1:z_2:z_3=\lambda {z'_1}^2:z'_1(z'_2-d
 z'_1):z'_3(z'_2-d z'_1), \label{eq:ex1'}
 \eea
 \be
 Jac(\Phi_2)=2 \lambda z_1 {z_2}^2,\qquad
 Jac(\Phi_2^{(-1)})=2\lambda {z'_1}^2(z'_2-d z'_1), \label{eq:jac1}
 \ee
 \bea
 O_1 & = & (1,0,0),\quad O_2^*=O_3^*=(0,0,1), \nonumber \\
  J_1 & : & (z_1=0), \quad J_2=J_3:(z_2=0),
 \label{eq:o1}
 \eea
 \bea
 O_1^{(-1)} & = &  (0,1,0) ,\quad {O_2^{(-1)}}^*={O_3^{(-1)}}^*=(0,0,1), \nonumber \\
 J_1^{(-1)} & : & (z'_2-d z'_1=0),\quad J_{2}^{(-1)}=J_{3}^{(-1)} :(z'_1=0). \label{eq:o1'}
\eea

\item{\bf the FE (90) from \citep{rerikh-92} $F(w+1)=\frac{3 F(w)-F(w-1)+
F(w)F(w-1)}{1+F(w)}$}.\label{ex2}

Omitting changes (see previous example) we have
 \be
 \Phi_2: z'_1:z'_2:z'_3  =  z_2(z_2+ z_3):3 z_2 z_3- z_1 z_3 + z_1 z_2:
  z_3(z_2+z_3),\label{eq:ex2}
\ee
\be
\Phi_2^{(-1)}: z_1:z_2:z_3  =  3 z'_1 z'_3- z'_2 z'_3 -z'_1 z'_2: z'_1(
z'_3-z'_1):z'_3(z'_3-z'_1) ,\label{eq:ex2'}
\ee
\be
 {\rm Jac}(\Phi_2)= 2(z_3-z_2)(z_2+z_3)^2,~
  {\rm Jac}(\Phi_2^{(-1)})=2 (z'_1+z'_3)(z'_3-z'_1)^2,
  \label{eq:jac2}
\ee
\bea
  O_1 & = & (-3/2,-1,1),\quad O_{2}^* = O_{3}^* = (1,0,0), \nonumber \\
J_1 & : & (z_3-z_2=0),\quad J_{2} = J_{3} : (z_3+z_2=0), \label{eq:o2} \\
O_1^{(-1)} & = & (1,3/2,1),\quad O_{2}^{(-1)*} = O_{3}^{(-1)*} = (0,1,0), \nonumber \\
  J_1^{(-1)} & : & ( z'_1+ z'_3=0),\quad J_{2}^{(-1)} = J_{3}^{(-1)} :(z'_3-z'_1=0) \label{eq:o2'}
\eea

 \item {\bf Mapping of the paper. }                                   \label{ex3}
\bea
\Phi_2 & : & z_1':z_2':z_3' = [z_1z_3-p_2z_1^2+\frac{(q_2+q_3)}{2}z_1z_2+
\frac{(q_2-q_3)^2}{12p_2}z_2^2]:\nonumber\\
& & -z_2[z_3+2p_2z_1+\frac{q_2+q_3}{2}z_2]:z_3^2+\frac{3}{2}(q_2+q_3)z_2z_3\nonumber\\
& & -p_2^2z_1^2+p_2(q_2+q_3)z_1z_2+\frac{1}{12}(5q_2^2+
14q_2q_3+5q_3^2)z_2^2,  \label{eq:c15}\\
\Phi_2^{-1} & = & \Lambda \comp \Phi_2 \comp \Lambda, \quad \Lambda=diag(-1,1,1). \label{eq:c16}
\eea
The mapping \rf{eq:c15} follows from generic quadratic map \rf{eq:quad} if we suppose $B_1= B\Lambda,
\quad \Lambda=diag(-1,1,1)$ and
\be
  B = \left(
 \begin{array}{ccc}
   p_1 & q_1 & r_1 \\
   p_2 & q_2 & r_2 \\
   p_3 & q_3 & r_3
 \end{array}
 \right),         \label{eq:c11}
\ee
where
\be
p_1=-2p_2,\quad p_3=p_2,\quad q_1=(q_2+q_3)/2, \quad r_1=r_2=r_3=1. \label{eq:c11a}
\ee

 In accordance with formulaes \rf{eq:quad}, \rf{eq:j} and \rf{eq:quadj}we have  three principal
 lines $J_i, J^{(-1)}_i$ and three $ F$-points $O_i, O^{(-1)}_i,~~O_i = (j_j=0)\cap
(j_k = 0),~~O^{(-1)}_i = (j^{(-1)}_j=0)\cap(j^{(-1)}_k = 0), ~i\neq j\neq k,~~
i,j,k\in(1,2,3)$:
\bea
  J_i &:& (j_i=- p_i z_1 + q_i z_2 + r_iz_3=0), \nonumber\\
  J^{(-1)}_i & : & (j^{(-1)}_i=
p_i z_1 + q_i z_2 + r_i z_3=0),  \label{eq:c12}\\
  O_i &=& \left\{ q_j - q_k, p_j - p_k, p_k q_j - p_j q_k \right\},
  ~~ ~~ O^{(-1)}_i = \Lambda O_i,\label{eq:c13}\\
O_1 &=&(\frac{1}{p_2},0,1),~O_2 =\left(-\frac{q_3-q_2}{p_2(5q_3+q_2)},-\frac{
6}{5q_3+q_2},1\right),\nonumber\\
& &O_3=\left(\frac{q_3-q_2}{p_2(5q_2+q_3)},-\frac{6}{5q_2+q_3},1\right).
\label{eq:c17}
\eea

\item {\bf FE $F(w+1)= \frac{4+2 F(w) F(w-1)+ F(w-1)-14 F^2(w)-4 F(w-1) F^2(w)}
{1-2 F(w)-2 F(w-1)-4 F^2(w)}$ of paper \citep{rerikh-95b}} \label{ex4}

Let us consider the cubic birational mapping  $\Phi_3: \quad {\rm\bf
CP^2\mapsto CP^2}$ associated with the above functional equation  from
\citep{rerikh-95b} (see eq. 23 on p. 67 and  eq. 30 on p. 68)
 \bea
\Phi_3 & : & z'_1:z'_2:z'_3 = z_2 (z^2_3-2 z_1 z_3-2 z_2 z_3-4 z^2_2):
 (4 z^3_3+z_1 z^2_3+\nonumber\\
 & ~ & 2 z_1 z_2 z_3 - 14 z^2_2 z_3 - 4 z_1 z^2_2):
 z_3 (z^2_3-2 z_1 z_3-2 z_2 z_3-4 z^2_2), \label{eq:c26}\\
 \Phi^{(-1)}_3 & : & z_1:z_2:z_3 = -(4 z'^3_3-z'_2 z'^2_3+2 z'_1 z'_2 z'_3
  - 14 z'^2_1 z'_3 + 4 z'_2
 z'^2_1): \nonumber \\
 & ~  & z'_1 (z^2_3+ 2 z'_1 z'_3+2 z'_2 z'_3-4 z'^2_1): \nonumber\\
 & ~ & z'_3 (z^2_3+2 z'_1 z'_3+2 z'_2 z'_3-4
 z'^2_1), \label{eq:c26a}
 \eea
where $y_1=F(w-1), y_1'= y_2=F(w), y_2'=F(w+1)~~ \mbox{and}~~ y_i=\frac{z_i}{z_3}$.

 These maps (char=\{3;2,1,1,1,1\}) have the following indeterminacy points and principal curves:
 \bea
O_1 &=& (1, 0, 0),~~O_2= (-1/2, 1/2, 1),~~O_3=(-5/2, -3/2, 1),\nonumber \\
 O_4 &=& (1/2,-1/2, 1),~~O_5= (-11/2, 3/2, 1), \label{eq:c27}\\
 O^{(-1)}_1 &=& (0,1,0),~~O^{(-1)}_2 = (3/2, 5/2, 1),~~O^{(-1)}_3
 =(-1/2, 1/2, 1), \nonumber \\
 O^{(-1)}_4 &=& (-3/2, 11/2, 1),~~O^{(-1)}_5=(1/2, -1/2, 1), \label{eq:c27a}\\
 J_1 &:& (z^2_3-2 z_1 z_3-2 z_2 z_3-4 z^2_2=0),~~J_2: (2 z_2-3 z_3=0), \nonumber \\
J_3 &:& (z_3+2 z_2=0),~ J_4 : (2 z_2+3 z_3=0),\nonumber\\
J_5 &:& (z_3-2z_2=0), \label{eq:c28} \\
J^{(-1)}_1 &:& (z'^2_3+2 z'_1 z'_3+2 z'_2 z'_3-4 z'^2_1=0),~~ J^{(-1)}_2
 : (2 z_1-z_3=0),\nonumber \\
  J^{(-1)}_3 &:& (2 z_1+ 3 z_3=0),\quad J^{(-1)}_4 : (2 z_1+z_3=0),\nonumber\\
  J^{(-1)}_5 & : & (2 z_1-3 z_3=0).\label{eq:c28a}
 \eea
 It is not difficult to obtain $i_{\alpha,\beta}$:
   \be
    i_{\alpha,\beta} =\left(
\begin{array}{ccccc}
              1 & 1 & 1 & 1 & 1\\
              1 & 0 & 0 & 0 & 1\\
              1 & 0 & 0 & 1 & 0\\
              1 & 0 & 1 & 0 & 0\\
              1 & 1 & 0 & 0 & 0
              \end{array}
             \right).
                          \label{eq:c28a'}
\ee

\end{enumerate}

\section{ Appendix B. Examples to Sections \ref{decomp}  and  \ref{dynamics.diff.eq}} \label{examples.a}

Below we return again to Examples \ref{ex1}-\ref{ex4} (see Appendix A ) for illustration
 of Sections \ref{decomp}  and
\ref{dynamics.diff.eq}. Firstly, we give the set ${\bf\rm O}^{(int)}$
derived with the help of the decomposition procedure, Section \ref{decomp} \rf{eq:f19},
\rf{eq:f19a}, and then give the equations of the dynamics \rf{eq:f22},
\rf{eq:f23} of maps \ref{ex1}-\ref{ex4}, and finish with a set of
difference equation for the Arnold complexity and its solution.

\begin{enumerate}

\item \label{ex1a}
This map \rf{eq:ex1} has two infinitely near points $O^*_2$ and $O^*_3 $.
($O^*_2$ is merged with $O^*_3$ in the direction $J_1$
The point $O_1$ belongs ${\bf\rm O}^{(inf)}:\quad
\Phi^{(-k)}_2(O_1)=(\lambda^k,-d \frac{\lambda^k-1}{\lambda-1}, 0),$ but,
say, the point $O^*_3 \in {\bf\rm O}^{(int)}:\quad O^*_3=O'^*_3$. We can
ascertain that the point $O^*_2$ does not belong to ${\bf\rm O}^{(int)}$:
 if we consider a map
$\Phi_{2,\epsilon}$ being a small deformation of initial map \rf{eq:ex1}
and having three different indeterminacy points ($ O_1=(1, 0, 0),~~
O_2=(0,\epsilon/\lambda, 1),~~O_3=(0, 0, 1)$)
 \be
  \Phi_{2,\epsilon}  :~~  z'_1:z'_2:z'_3=z_2 z_1:z_2(\lambda z_2-\epsilon z_3
  +d z_1):z_1 z_3, \label{eq:ex01a}
  \ee
Since it is difficult to point out a general method of constructing a map
coinciding with a given map at a small parameter $\epsilon=0$ and having only
ordinary indeterminacy points, we change the given map by a birationally
equivalent map with ordinary indeterminacy points. The map $\Phi_2:
z\mapsto z',\quad z, z' \in {\bf\rm CP}^2$ \rf{eq:ex1} is birationally
equivalent to the map $\Phi_3: u\mapsto u',~u'\in {\bf\rm CP}^2, \quad
\Phi_3=\Psi^{(-1)}_2\comp\Phi_2\comp\Psi_2$
 \bea
  \Phi_3 & : &~ u'_1:u'_2:u'_3 = [(\lambda-d)u_1+d u_2][-u_3((\lambda-d)u_1+
 d u_2) + \nonumber \\
&   &~ (u_2-u_1)(u_3+u_2)]: [(1+d)u_2- (1+d-\lambda)u_1][-u_3((\lambda \nonumber \\
&   &~ -d)u_1 +  d u_2)+  (u_2-u_1)(u_3+u_2)]:\nonumber \\
&   &~ u_3[(\lambda-d)u_1+d u_2][(1+d)u_2-(1+d- \lambda)u_1], \label{eq:ex11a}\\
\Phi^{(-1)}_3 & : &~ u_1:u_2:u_3 = [u'_1(d+1)-d u'_2][2 u'_1 u'_3+u'_1u'_2-u'_2u'_3]:   \nonumber\\
&   &~ [u'_1(d+1-\lambda)- (d-\lambda)u'_2]  [2 u'_1u'_3+u'_1u'_2-u'_2 u'_3]:  \nonumber \\
&   &~   u'_3(u'_2-u'_1)[u'_1(d+1-\lambda)-(d-\lambda)u'_2],\label{eq:ex11a'}
 \eea
where the map $\Psi_2:\quad u\mapsto z$ is chosen so that two
indeterminacy points of the map $\Psi^{(-1)}_2$ may coincide with  the points
$O_1, ~O^*_2$ from \rf{eq:o1}, but  a direction of the second principal
curve (line) for $\Psi^{(-1)}_2$ through the point $O^*_2$ does not coincide
with $J_1$ from \rf{eq:o1}:
 \be
  \Psi_2:~u\mapsto z ~z_1:z_2:z_3=(u_2 u_3-u_1 u_3):u_1 u_3:
  (-u_1 u_3+u_1 u_2). \label{eq:ex12a}
  \ee
The maps $\Phi_3,~\Phi^{(-1)}_3 $ have  $ {\rm char}
=\{2,1,1,1,1\}$  and the following indeterminacy points:
 $$ O_1=(0,0,1),\qquad O_2=(1,1,0),\qquad  O_3=(1,0,0),$$
$$ O_4=(-\frac{d}{d-\lambda},-1,1),\qquad O_5=(-\frac{2(1+d)}{1-\lambda+d},-2, 1),$$
$$ O^{(-1)}_1 =(0,0,1),\qquad O^{(-1)}_2=(1,0,0),\qquad O^{(-1)}_3=(0,1,0),$$
$$O^{(-1)}_4=(-\frac{d-\lambda-1}{d-\lambda+1},-\frac{d-\lambda-1}
{d-\lambda},1),\qquad O^{(-1)}_5=(-1,-1,1).$$

 Making the decomposition we
obtain
 $$
 O_1, O_3  \in  {\bf\rm O}^{(int)}, ~ \Phi^{(-1)}_3 O_2=O_2,~
  O_2\in{\bf\rm O}^{(cycle)}, ~ O_4, O_5\in {\bf\rm O}^{(inf)},
 $$
$$
  O^{(-1)}_1, O^{(-1)}_2  \in  {\bf\rm O}^{(-1)(int)}, ~
  O^{(-1)}_3, O^{(-1)}_4 \in {\bf\rm O}^{(-1)(inf)},~ O^{(-1)}_5\in
   {\bf\rm O}^{(-1)(cycle)}.
$$
 Following Theorem \ref{dynamics} we have ( $ i_{\alpha \beta}$ is the same one
 as in example \ref{ex4} \rf{eq:c28a'})
 \bea
 d(k) & = & 3 d(k-1)-2 \gamma_1(k-1)-\gamma_3(k-1), \label{eq:ex14a}\\
\gamma_1(k) & = & 2 d(k-1)-\gamma_1(k-1)-\gamma_3(k-1),\\
 \gamma_3(k) & = & d(k-1)-\gamma_1(k-1).
 \eea
 Following Theorem \ref{dynamics} we obtain
 \be
 d(k+2)-2 d(k+1)+d(k)=0,\qquad d(k)=2 k+1.\label{eq:ex15a}
 \ee

 \item \label{ex2a}

The map $\Phi_2$ \rf{eq:ex2} is birationally equivalent to  the map
$\Phi^*_2=\Psi^{(-1)}_2\comp\Phi_2\comp\Psi_2$ with ordinary indeterminacy
points ( conditions for the choice of $\Psi_2$:
$O^{(-1)}_1(\Psi^{(-1)}_2)=O_{2,3}(\Phi_2),~
O^{(-1)}_2(\Psi^{(-1)}_2)=O^{(-1)}_{2,3}(\Phi^{(-1)}_2)$, $J^{(-1)}_1,
J_1~\mbox{for}~\Phi_2$ from \rf{eq:o2}, \rf{eq:o2'} must not coincide with
$J_i $ for $\Psi_2$)
 \bea
\Phi^*_2  : \quad u'_1:u'_2:u'_3 & = & u_1(3 u_2-u_3+u_1):u_3(3
u_2-u_3+u_1):\nonumber\\
&  & u_2(u_1+u_3), \label{eq:ex2a}\\
\Phi^{*(-1)}_2  : \quad u_1:u_2:u_3 & = & u'_1(u'_2+u'_1-3
u'_3):u'_3(u'_1-u'_2):\nonumber\\
&  & u'_2(u'_2+u'_1-3 u'_3),\label{eq:ex2a'}
 \eea
 but the map $\Psi_2$ is
 $
 \Psi_2:\quad z_1:z_2:z_3=u_1u_3:u_1u_2:u_2u_3.
$

We have the indeterminacy points for $\Phi^*_2$:~
  $O_1  =  (-1,\frac{2}{3},1),~ O_2=(1,0,1), $ ~~ $ O_3=(0,1,0), ~
  O^{(-1)}_1 = (\frac{3}{2}, \frac{3}{2}, 1),~  O^{(-1)}_2=(-1, 1, 0), ~
  O^{(-1)}_3=(0, 0, 1)$.

  The decomposition of the sets ${\bf\rm O, O}^{(-1)}$ gives:~
  $  O_1, O_2 \in{\bf\rm O}^{(inf)},\quad O^{(-1)}_1,$ ~ $ O^{(-1)}_2  \in {\bf\rm O}^{(-1)(inf)},$ \quad
  $ \Phi^{*(-1)}_2(O_3)=O^{(-1)}_3,~ O_3\in{\bf\rm O}^{(int)},~ O^{(-1)}_3\in{\bf\rm O}^{(-1)(int)}.$

  We have for the Arnold complexity
  $  d(k)  =  2 d(k-1)-\gamma_3(k-2),  \gamma_3(k)  =  d(k-1)$,\quad
 $
  d(k+3)-2 d(k+2)+d(k)=0,\quad d(k)=-1 +(\lambda^{k+3}+
  (-1)^k\lambda^{-(k+3)})/\sqrt 5,$
  where $\lambda=\frac{\sqrt 5+1}{2}$.
\item \label{ex3a}
The decomposition of the sets  ${\bf\rm O, O}^{(-1)}$ of indeterminacy points $O_\alpha,
 O^{(-1)}_\beta$ \rf{eq:c13}, \rf{eq:c17} of the mappings ${\bf\rm \Phi_2}$ \rf{eq:c15} and
${\bf\rm \Phi^{(-1)}_2}$ \rf{eq:c16} gives
\bea
\Phi_2^{-1}(O_1) & = &\{\infty, 0\},~~~~~~~~~~~~~~~~~~~
 \qquad \Phi_2^{-2}(O_1) = O^{(-1)}_1, \nonumber \\
\Phi_2^{-1}(O_2) & = &\{\frac{q_3-q_2}{3p_2(q_3+q_2)},-\frac{2}{q_3+q_2}\},
~\qquad \Phi_2^{-2}(O_2) = O^{(-1)}_3, \nonumber \\
\Phi_2^{-1}(O_3) & = &\{-\frac{q_3-q_2}{3p_2(q_3+q_2)},-\frac{2}{q_3+q_2}\}, \quad
\Phi_2^{-2}(O_3) = O^{(-1)}_2,  \label{eq:ex3a1}
\eea
and, consequently,
${\bf\rm O}\equiv{\bf\rm O}^{(int)},~ {\bf\rm O}^{(-1)}\equiv{\bf\rm O}^{(-1)(int)},~
m_j=2 ~ \forall ~ j\in(1,2,3),$ $ \alpha_j=(1,2,3),~ \beta_j=(1,3,2).$

Due to Theorem \rf{dynamics} and Remark \rf{equal.some.m_i} we have from \rf{eq:f21}-
\rf{eq:sec.eq}
\bea
d(k)=2d(k-1)-S(k-3),  \label{eq:ex3a3} \\
S(k)=3d(k-1)-2S(k-3),  \label{eq:ex3a4}
\eea
where $S(k)= \sum_{\alpha=1}^{3}\gamma_\alpha(k)$.
In correspondence with Theorem \ref{dynamics} we obtain the difference equation
for $d(k)$
\be
d(k+4)-2d(k+3)+2d(k+1)-d(k)=0   \label{eq:ex3a5}
\ee
and its general solution
\be
d(k)=\frac{3}{4}k^2-\frac{1}{8}(-1)^k +\frac{9}{8}. \label{eq:ex3a6}
\ee
\item \label{ex4a}
The decomposition of the sets  ${\bf\rm O, O}^{(-1)}$ of indeterminacy points $O_\alpha,
 O^{(-1)}_\beta$ \rf{eq:c27}, \rf{eq:c27a} of the mappings ${\bf\rm \Phi_3}$ \rf{eq:c26} and
${\bf\rm \Phi^{(-1)}_3}$ \rf{eq:c26a} gives:
\bea
O_2 & = & O^{(-1)}_3, \quad  \alpha_1=2, \beta_1=3, \quad m_1=0, \nonumber \\
O_4 & = & O^{(-1)}_5, \quad \alpha_2=4, \beta_2=5, \quad m_2=0, \nonumber \\
\Phi^{(-1)}_3(O_1) & = & O^{(-1)}_1, \quad \alpha_3=1, \beta_3= 1, \quad m_3=1, \label{eq:ex3a7}\\
{\bf\rm O}^{(int)} & = & \{O_2, O_4, O_1\}, \quad
 {\bf\rm O}^{(-1)(int)}=\{O^{(-1)}_3, O^{(-1)}_5, O^{(-1)}_1\}, \nonumber \\
{\bf\rm O}^{(inf)} & = & \{O_3, O_5\} ,\quad {\bf\rm O}^{(-1)(inf)}= \{ O^{(-1)}_2, O^{(-1)}_4\}. \nonumber
\eea
Following Theorem \ref{dynamics} (see \rf{eq:sec.eq} and\rf{eq:lambda}) we have
$ Det(\Lambda) = (\lambda-1)^2(\lambda^3-\lambda^2-\lambda-1)$ and
$d(k)=c_{00}+c_{01}k+\sum_{i=1}^{i=3}c_{1i}\lambda^k_i$, where the coefficients $c_{0i}$
 and $c_{1i}$ are defined
in terms of $\lambda_i$--the roots of the cubic equation $\lambda^3-\lambda^2-\lambda-1=0$.
\end{enumerate}

\newpage

\begin{thebibliography}{00}
\bibitem[{Anosov et~al.(1988)Anosov, Bronshtein, Aranson, and
  Grines}]{DS-1-88b}
Anosov, D.~V., Bronshtein, I.~U., Aranson, S.~K., Grines, V.~Z., 1988. Smooth
  Dynamical Systems. Eds. D.~V.~Anosov and V.I.~Arnold, Dynamical Systems, Vol.
  1, {E}ncyclopaedia Math. Sciences. Springer, Berlin, pp. 151--242.

\bibitem[{Arnold(1984)}]{arnold-88c}
Arnold, V.~I., 1984. Ordinary differential equations. Springer-Verlag, Berlin.

\bibitem[{Arnold(1988)}]{arnold-88}
Arnold, V.~I., 1988. Geometrical Methods in the Theory of Ordinary Differential
  Equations. Springer-Verlag, New York.

\bibitem[{Arnold(1990{\natexlab{a}})}]{arnold-90b}
Arnold, V.~I., 1990{\natexlab{a}}. Dynamics of intersections. In: Rabinowitz,
  P., Zehnder, E. (Eds.), Proceedings of a Conference in Honour of J.~Moser.
  Academic Press, New York, pp. 77--84.

\bibitem[{Arnold(1990{\natexlab{b}})}]{arnold-90a}
Arnold, V.~I., 1990{\natexlab{b}}. Dynamics of the complexity of intersections.
  Bol. Soc. Bras. Mat. 21, 1--10.

\bibitem[{Arnold and Il'yashenko(1988)}]{DS-1-88a}
Arnold, V.~I., Il'yashenko, Y.~S., 1988. Ordinary Differential Equations. Eds.
  D.~V.~Anosov and V.I.~Arnold, Dynamical Systems, Vol.1, {E}ncyclopaedia Math.
  Sciences. Springer, Berlin, pp. 7--148.

\bibitem[{Baker(1966)}]{baker-66}
Baker, A., 1966. Linear forms in logarithms of algebraic numbers. Matematika
  13, 204--216.

\bibitem[{Baker(1967{\natexlab{a}})}]{baker-67a}
Baker, A., 1967{\natexlab{a}}. Linear forms in logarithms of algebraic numbers,
  ii. Matematika 14, 102--107.

\bibitem[{Baker(1967{\natexlab{b}})}]{baker-67b}
Baker, A., 1967{\natexlab{b}}. Linear forms in logarithms of algebraic numbers,
  iii. Matematika 14, 220--228.

\bibitem[{Baker(1968)}]{baker-68}
Baker, A., 1968. Linear forms in logarithms of algebraic numbers, iv.
  Matematika 15, 204--216.

\bibitem[{Baker(1971)}]{baker-71}
Baker, A., 1971. Effective methods in the theory of numbers. In: Proceedings of
  the International Congress of Mathematicians, Nice, September 1970. Vol.~1.
  Gauthier-Villars, 55, quai des Grands-Augustins, Paris $\mbox{6}^e$, pp.
  19--26.

\bibitem[{Baker(1990)}]{baker-90}
Baker, A., 1990. Transcendental Number Theory. Cambridge University Press,
  Cambridge.

\bibitem[{Baker and W\"ustholz(1993)}]{baker.wustholz-93}
Baker, A., W\"ustholz, G., 1993. Logarithmic forms and group varieties. J.
  reine anngew. Math. 442, 19--62.

\bibitem[{Coble(1961)}]{coble-61}
Coble, A.~B., 1961. Algebraic Geometry and Theta Functions. AMS, Providence, R.
  I., {A}MS Colloquium Publications, vol. X.

\bibitem[{Feldman(1968)}]{feldman-68-rus}
Feldman, N.~I., 1968. Improvement of evaluation of the linear form of
  logarithms of algebraic numbers (in russian). Mat. Sbornik 77 (119) 3,
  423--436.

\bibitem[{Feldman(1982)}]{feldman-82-rus}
Feldman, N.~I., 1982. Hilbert's Seventh Problem. Moscow University Press,
  Moscow.

\bibitem[{Gel'fond(1971)}]{gel'fond-67}
Gel'fond, A.~O., 1971. Calculus of Finite Differences. Hindustan Publishing
  Corporation, Delhi, authorized English translation of the third Russian
  edition.

\bibitem[{Griffiths and Harris(1978)}]{griffiths.harris-78}
Griffiths, P., Harris, J., 1978. Principles of Algebraic Geometry. John Wiley
  \& Sons, New York.

\bibitem[{Hudson(1927)}]{hudson-27}
Hudson, H., 1927. Cremona Transformations in Plane and Space. Cambridge
  University Press, Cambridge.

\bibitem[{Iskovskikh and Reid(1991)}]{iskovskikh.reid-91}
Iskovskikh, V.~A., Reid, M., 1991. Foreword to Hudson's book ``Cremona
  transformations''. Cambridge University Press, Cambridge, {U}npublished.

\bibitem[{Moser(1960)}]{moser-60}
Moser, J., 1960. On the integrability of area-preserving cremona mappings near
  an elliptic fixed point. Bol. Soc. Mat. Mexicana, 176--180.

\bibitem[{Moser(1994)}]{moser-94}
Moser, J., 1994. On quadratic symplectic mappings. Mathematische Zeitschrift
  216, 417--430.

\bibitem[{Rerikh(1992)}]{rerikh-92}
Rerikh, K.~V., 1992. Cremona transformation and general solution of one
  dynamical system of the static model. Physica D 57, 337--354.

\bibitem[{Rerikh(1995{\natexlab{a}})}]{rerikh-95a}
Rerikh, K.~V., 1995{\natexlab{a}}. Non-algebraic integrability of one
  reversible cremona dynamical system. the poincare (1.1) resonance and the
  birkhoff-moser analytical invariants. In: Proc. of Inter. Workshop "Finite
  dimensional integrable systems". JINR, Dubna, pp. 171--180.

\bibitem[{Rerikh(1995{\natexlab{b}})}]{rerikh-95b}
Rerikh, K.~V., 1995{\natexlab{b}}. Non-algebraic integrability of the chew-low
  reversible dynamical system of the cremona type and the relation with the 7th
  hilbert problem (non-resonant case). Physica D 82, 60--78.

\bibitem[{Rerikh(1997)}]{rerikh-97}
Rerikh, K.~V., 1997. Algebraic addition concerning the siegel theorem on the
  linearization of a holomorphic mapping. Math.Z. 224, 445--448.

\bibitem[{Rerikh(1998{\natexlab{a}})}]{rerikh-98a}
Rerikh, K.~V., 1998{\natexlab{a}}. Algebraic-geometry approach to integrability
  of birational plane mappings. {I}ntegrable birational quadratic reversible
  mappings. {I}. J. of Geometry and Physics 24, 265--290.

\bibitem[{Rerikh(1998{\natexlab{b}})}]{rerikh-98b}
Rerikh, K.~V., 1998{\natexlab{b}}. Non-algebraic integrability of one
  reversible dynamical system of the cremona type. J. of Math. Phys. 39,
  2821--2832.

\bibitem[{Shafarevich(1977)}]{shafarevich-77}
Shafarevich, I.~R., 1977. Basic Algebraic Geometry. Springer, Berlin.

\bibitem[{Snyder et~al.(1970)Snyder, Coble, Emch, Lefschetz, Sharpe, and
  Sisam}]{topics-70}
Snyder, V., Coble, A.~B., Emch, A., Lefschetz, S., Sharpe, F.~R., Sisam, C.~H.,
  1970. Selected Topics in Algebraic Geometry. Chelsea, Bronx, N.Y., "A reprint
  in one volume, with the correction of errata, of two works published under
  the title, Selected topics in algebraic geometry, Bulletin of the National
  Research Council, number 63 (Washington, 1928); Bulletin of the National
  Research Council, number 96 (Washington, 1934). The original list of books
  has been supplemented by inclusion of later editions and printings, and
  collected works.".

\bibitem[{Veselov(1989)}]{veselov-89}
Veselov, A.~P., 1989. Cremona group and dynamical systems. Mat. zametki 45, 3,
  118--120.

\bibitem[{Veselov(1991)}]{veselov-91}
Veselov, A.~P., 1991. Integrable mappings. Russian Math. Surveys 46, no. 5,
  1--51.

\bibitem[{Veselov(1992)}]{veselov-92}
Veselov, A.~P., 1992. Growth and integrability in the dynamics of mappings.
  Commun. Math. Phys. 145, 181--193.

\end{thebibliography}

\end{document}